\documentclass{amsart}
\usepackage{bbm}
\usepackage{color}
\usepackage{amscd}
\usepackage{nicefrac}
\usepackage{graphicx}
\usepackage{hyperref}
\usepackage{amssymb}
\usepackage{euscript}
\usepackage{enumitem}
\usepackage[cmtip,all]{xy}
\usepackage[export]{adjustbox}
\usepackage{slashed}

\numberwithin{equation}{section}
\renewcommand{\subsection}[1]{\hspace{-\parindent}\refstepcounter{subsection}{\bf (\thesubsection) #1.}}
\renewcommand{\thesubsection}{\arabic{section}\alph{subsection}}

\allowdisplaybreaks

\theoremstyle{plain}
\newtheorem{thm}{Theorem}[section]

\newtheorem{theorem}[thm]{Theorem}

\newtheorem{corollary}[thm]{Corollary}
\newtheorem{assumptions}[thm]{Assumptions}

\newtheorem{lemma}[thm]{Lemma}

\newtheorem{remark}[thm]{Remark}

\newtheorem{non-example}[thm]{Non-example}

\newtheorem*{claim*}{Claim} 
\newtheorem*{lemma*}{Lemma}
\newtheorem*{theorem*}{Theorem}
\newtheorem*{conjecture*}{Conjecture}

\newcommand{\bC}{{\mathbb C}}

\newcommand{\bF}{{\mathbb F}}

\newcommand{\bR}{{\mathbb R}}

\newcommand{\bZ}{{\mathbb Z}}

\newcommand{\frakP}{\mathfrak{P}}

\newcommand{\iso}{\cong}
\newcommand{\htp}{\simeq}

\newcommand{\rel}{\operatorname{rel}}
\newcommand{\dq}{\slashed{\partial}_q}
\newcommand{\nablaq}{\slashed{\nabla}_q}

\begin{document}

\title[$p$-curvature]{$P$-curvature from a TQFT viewpoint}
\author{Paul Seidel}
\maketitle

\begin{abstract}
This paper is motivated by the recent discovery, in symplectic topology, of a relation (partly proved, partly conjectural) between quantum connection and quantum Steenrod operations. We introduce an abstract TQFT framework, based on certain features of the definition of deformed symplectic cohomology,  in which such a relation holds generally. 
\end{abstract}

\section{Introduction}
$P$-curvature is an invariant of algebraic connections (first order linear differential operators) over fields of characteristic $p>0$. It is not related to the ordinary notion of curvature, but it plays a similar role: vanishing of the $p$-curvature ensures the existence of ``sufficiently many'' covariantly constant local sections \cite[Theorem 5.1]{katz70}. More recently, $p$-curvature has been considered for the quantum connection from Gromov-Witten theory. The guiding conjecture, introduced by J.-H.~Lee with input from Etingof, can be stated informally as follows:
\begin{quote}
{\em Take the quantum connection, which differentiates in the direction of some degree $2$ cohomology class. Its $p$-curvature agrees with the quantum Steenrod endomorphism associated to that class.}
\end{quote}
The groundbreaking paper \cite{lee23b} proves it in an important instance involving equivariant quantum cohomology, namely the cotangent bundle of a flag variety \cite[Corollary 5.2]{lee23b}. This was motivated by the appearance of that particular instance of $p$-curvature in representation theory \cite{etingof-varchenko24}. For Fano varieties, a weaker version of the conjecture (agreement up to a nilpotent error) was used in \cite{chen24b} to prove that the quantum connection has a singularity of exponential type at $\infty$; in the same situation, the argument was improved to a full proof of the conjecture in \cite[Proposition 1.7]{pomerleano-seidel26}. A proof for Calabi-Yau threefolds was obtained in \cite[Theorem 1.5]{bai-lee-pomerleano26}. The paper \cite{rezchikov26} proves it for Calabi-Yau's of any dimension, under a number of assumptions on the Fukaya category, and assuming that the prime $p$ is large. Conjecturally, there are potential generalizations to a relative situation involving a smooth ample divisor \cite[Remark 1.15]{pomerleano-seidel26}, and to the fixed point Floer cohomology of $p$-th iterates of symplectic automorphisms.

Particularly after looking at the last-mentioned generalizations, it is natural to ask if the relation between quantum connection and Steenrod-type operations can be viewed as a general property of a class of two-dimensional ``closed string TQFTs'', which means chain complexes which carry operations parametrized by chains on suitable moduli spaces of Riemann surfaces without boundary. (A negative answer could mean that the conjecture only holds in the presence of a much richer ``open string'' categorical context, as in \cite{rezchikov26}; or, that it is a truly geometric insight depending on a sophisticated analysis of pseudo-holomorphic curves, as in \cite{bai-lee-pomerleano26}.) Because the conjecture addresses a number of different contexts, there is not necessarily a unique answer. This paper argues in favor of the ``closed string TQFT'' view. We introduce an abstract framework of  operations parametrized by cylinders with marked points. Within that framework, one can define versions of the quantum connection and quantum Steenrod operations, and their relation is indeed as predicted by the guiding conjecture. The general result is Theorem \ref{th:p-curvature-conjecture}, which is then transformed into two versions closer to the traditional statement (Corollaries \ref{th:p-curvature-module} and \ref{th:p-curvature-module-2}).

The approach here is motivated by the approach to quantum cohomology relative to a smooth ample divisor, or more technically deformed symplectic cohomology, in \cite{pomerleano-seidel24}. That approach involves counting solutions of Floer-type equations on cylinders with marked points, where the marked points dictate the intersection points with an ample divisor in the target manifold. Cylinders with $m$ marked points are counted with $q^m$, where $q$ is the quantum or Novikov variable. The marked points are not numbered, and can collide freely, so that the parameter spaces are symmetric products of the cylinder. Our setup has the same formal structure, starting with a $q$-deformed differential. In fact, we will consider two versions of the theory.
\begin{itemize} \itemsep.5em
\item The first version corresponds to $S^1$-equivariant Floer cohomology (even though, for technical reasons, we replace $S^1$ by an infinite discrete subgroup); this is constructed in Section \ref{sec:tqft}. $S^1$-equivariant Floer cohomology is the standard setting for constructing the quantum connection, see e.g.\ \cite[Section 5.3f]{pomerleano-seidel23}; we will carry out an abstract version of that construction in Section \ref{sec:connection}, and then directly compute the $p$-curvature. 

\item The second version, in Section \ref{sec:tqft-2}, uses only a $\bZ/p$ group of symmetries, which is the standard context for quantum Steenrod operations. It is related to the $S^1$-equivariant one in a familiar way \cite{sen23}, and the natural endpoint of the argument is to translate the outcome of the $p$-curvature context accordingly. (It would have been possible to carry out the construction of the connection and the $p$-curvature computation entirely in $\bZ/p$-equivariant terms, see Remark \ref{th:shortcoming-2}; but that is not particularly natural when thinking of the quantum connection.)
\end{itemize}
The abstract nature of the argument allows us to formulate it in elementary terms, using only basic algebraic topology (singular chains and the chain level cross product) and homological algebra ($q$-filtration arguments).

We should discuss the relation between this paper and \cite{chen26}. There is a significant difference in the aims and underlying geometry: because \cite{chen26} focuses on applications to cyclic homology (similarly to \cite{rezchikov26}, but without the latter's use of cyclotomic structures); for that reason, it necessarily uses the Kontsevich-Soibelman operad of configuration spaces of points on the cylinder. Those have different topology from our symmetric products, which come from the Floer-theoretic constructions in \cite{pomerleano-seidel24}. Nevertheless, the structure of the two theories is similar in many respects. The appearance of two kinds of ``quantum Steenrod'' operations (for us, this is \eqref{eq:rq} and \eqref{eq:rhoq} in the $S^1$-equivariant case, as well as \eqref{eq:sq} and \eqref{eq:sigmaq} in the $\bZ/p$-equivariant case) mirrors \cite[Theorem 5.8]{chen26}. The relation between operations in the two flavors of the theory, Lemma \ref{th:compare-operations}, is parallel to \cite[Corollary 5.14]{chen26}. However, \cite{chen26} contains no result corresponding to our $p$-curvature computations. 

The reader may have noticed the weasel-wordy ``argue for'' above. There are no applications to symplectic topology in this paper. The natural target would be the quantum connection on closed Calabi-Yau manifolds, via an adaptation of the approach in \cite{pomerleano-seidel24}. However, checking whether the relevant Floer theory fits into our abstract framework, or into an equivalent version of it, would require a significant amount of additional work. This is in part a general issue with chain level operations (compare \cite{abouzaid-groman-varolgunes22}, which established the chain level $E_2$ structure of symplectic cohomology). A more specific problem is that, as written here, group actions on chain complexes are strict, something that is not true for the standard Floer complexes. However, there is a significant mitigating aspect: all our operations are constructed by simple cohomological obstruction theory arguments. Hence, one should be able to replace the framework introduced here by another one with the same cohomology level behaviour; this leaves space for compromises dictated by Floer theory practicalities, while preserving the same structure of proofs.

{\em Acknowledgments.} I am indebted to Dan Pomerleano for discussions that spawned the idea pursued in this paper, and to Zihong Chen for explaining aspects of \cite{chen26}.

{\em AI disclosure.} During the project, AI (ChatGPT) provided the formula \eqref{eq:chain-exterior-product-3}, derived from \cite{maclane}. After completion of an earlier version of the manuscript, AI review pointed out two significant mistakes: a conceptual one (Remark \ref{th:shortcoming-0}), which caused me to revert to a previous variant of the overall strategy, using only discrete groups on the chain level; and a technical one, whose fix (Lemma \ref{th:order-zero}) includes a suggestion by the AI.

\section{Equivariant cohomology}
Our discussion of equivariant cohomology uses 1950s style explicit chain level formulae. This keeps the discussion elementary, but has limitations which will become visible as we proceed (Remarks \ref{th:shortcoming-0} and \ref{th:shortcoming-1}). Presumably, $\infty$-categorical language, as in \cite{chen24b}, would do better.

\subsection{Equivariant cohomology}
Let $\Gamma$ be a (discrete) group. In this section, all complexes are ones of vector spaces over some field $\bF$. We consider the group ring $\bF \Gamma$, and its augmentation
\begin{equation} \label{eq:augmentation}
\epsilon: \bF \Gamma \longrightarrow \bF.
\end{equation}
Given a chain complex $V$ on which $\Gamma$ acts, we write $H^*_\Gamma(V)$ for equivariant cohomology. One can define it as 
\begin{equation} \label{eq:derived-hom}
H^*_\Gamma(V) = \mathit{DHom}_{\bF \Gamma}(\bF,V),
\end{equation}
the space of morphisms in the derived category (as a graded space, meaning that it includes the morphisms to all shifted versions). For $V = \bF$ with the trivial action, $H^*_\Gamma(\bF)$ is group cohomology with $\bF$-coefficients. Equivariant cohomology comes with canonical maps
\begin{equation} \label{eq:include-and-forget}
H^*(V^\Gamma) \longrightarrow H_\Gamma^*(V) \longrightarrow H^*(V),
\end{equation}
where $V^\Gamma$ is the fixed subcomplex. It also has an exterior product structure
\begin{equation} \label{eq:exterior-product}
\times: H^*_\Gamma(V) \otimes H^*_\Gamma(W) \longrightarrow H^*_\Gamma(V \otimes W),
\end{equation}
which is associative, (graded) commutative, and compatible with \eqref{eq:include-and-forget} in the obvious way.  In particular, $H^*_\Gamma(V)$ is always a module over $H^*_\Gamma(\bF)$; for that special case, we will usually omit the $\times$ from the notation. The origin of the exterior product can be explained in a way similar to \eqref{eq:tensor-product}, as a composition
\begin{equation} \label{eq:exterior-product-2}
\begin{aligned}
& \mathit{DHom}_{\bF \Gamma}(\bF,V) \otimes \mathit{DHom}_{\bF \Gamma}(\bF, W)
\longrightarrow \mathit{DHom}_{\bF \Gamma^2}(\bF, V \otimes W)
\\ & \qquad \qquad \qquad
\longrightarrow \mathit{DHom}_{\bF \Gamma}(\bF, V \otimes W),
\end{aligned}
\end{equation}
the last step being restriction to the diagonal subgroup.

We will also use a relative version $H^*_{\Gamma \rel}(V)$ of equivariant cohomology, which sits in a long exact sequence
\begin{equation} \label{eq:relative-les}
\cdots \rightarrow H^*_{\Gamma \rel}(V) \longrightarrow H^*_\Gamma(V) \longrightarrow H^*(V) \rightarrow \cdots
\end{equation}
One way to construct it is to think of $H^*(V) = D\mathit{Hom}_{\bF\Gamma}(\bF\Gamma,V)$. From this point of view, the forgetful map in \eqref{eq:include-and-forget} is given by composition with \eqref{eq:augmentation}. One would then consider the mapping cone, or equivalently up to shift and quasi-isomorphism, its kernel:
\begin{equation} \label{eq:relative-derived-hom}
H^*_{\Gamma \rel}(V) = D\mathit{Hom}_{\bF \Gamma}(\mathit{ker}(\epsilon)[1], V)
\end{equation}
(the chain level construction we will use later on is less theoretically satisfying, but simpler). The relative theory comes with versions of the exterior product,
\begin{equation} \label{eq:relative-product}
\begin{aligned}
& H^*_{\Gamma\rel}(V) \otimes H^*_{\Gamma}(W) \longrightarrow 
H^*_{\Gamma\rel}(V \otimes W), \\
& H^*_{\Gamma}(V) \otimes H^*_{\Gamma\rel}(W) \longrightarrow
H^*_{\Gamma\rel}(V \otimes W),
\end{aligned}
\end{equation}
with the same associativity and commutativity properties as in the absolute situation.

In the topological situation, where $\Gamma$ acts on a topological space $X$, our applications will concern the complex of singular chains $V = C_{-*}(X)$ (to match the general preference for cohomological conventions, the grading has been reversed). The resulting $H^*_\Gamma(V)$ does not reproduce equivariant homology or cohomology in the ordinary sense (a better parallel are the equivariant Chow groups in \cite{edidin-graham97}). However, if $X$ is a compact oriented $n$-manifold with boundary, and the $\Gamma$-action preserves the orientation, then Poincar{\'e} duality gives an isomorphism
\begin{equation} \label{eq:poincare-duality}
\xymatrix{
\ar[d]
H^*_\Gamma(C_{-*}(X)) \ar[r]^-{\iso} & H^*_\Gamma(C^{*+n}(X,\partial X)) \ar[r]^-{\iso} & H^{*+n}_\Gamma(X,\partial X) \ar[d]
\\
H_{-*}(X) \ar[rr]^{\iso} && H^{*+n}(X,\partial X)
}
\end{equation}
If $X$ and $Y$ are two spaces with $\Gamma$-actions, the standard cross product 
\begin{equation} \label{eq:cross-product}
C_{-*}(X) \otimes C_{-*}(Y) \longrightarrow C_{-*}(X \times Y) 
\end{equation}
and \eqref{eq:chain-exterior-product-2} induce a map
\begin{equation} \label{eq:topological-exterior-product}
\times: H^*_\Gamma(C_{-*}(X)) \otimes H^*_\Gamma(C_{-*}(Y)) \longrightarrow H^*_\Gamma(C_{-*}(X \times Y)).
\end{equation}

\subsection{Chain level constructions}
We want to fix (functorial) chain complexes underlying $H^*_\Gamma(V)$, denoted by $V^{h\Gamma}$. These should come with chain level realizations of \eqref{eq:include-and-forget} and \eqref{eq:exterior-product}, which are maps
\begin{align} 
\label{eq:chain-include-and-forget}
& V^{\Gamma} \longrightarrow V^{h\Gamma} \longrightarrow V, \\
\label{eq:chain-exterior-product}
& \times: V^{h\Gamma} \otimes W^{h\Gamma} \longrightarrow (V \otimes W)^{h\Gamma}.
\end{align}
It will be convenient for us to require that the chain level model has certain properties. First, the map \eqref{eq:chain-exterior-product} is strictly unital (with unit element in $\bF^{h\Gamma}$), strictly associative, and compatible with \eqref{eq:chain-include-and-forget}. Our second requirement is that the forgetful map \eqref{eq:chain-include-and-forget} is surjective. We can then define the relative version using
\begin{equation} \label{eq:kernel}
V^{h\Gamma \rel} \stackrel{\mathrm{def}}{=} \mathit{ker}(V^{h\Gamma} \rightarrow V).
\end{equation}
Because of the compatibility between product and forgetful maps, \eqref{eq:chain-exterior-product} restricts to maps
 \begin{equation} \label{eq:chain-relative-product}
 \begin{aligned}
\times: V^{h\Gamma \rel} \otimes W^{h\Gamma} \longrightarrow (V \otimes W)^{h\Gamma \rel}, \\
\times: V^{h\Gamma} \otimes W^{h\Gamma \rel} \longrightarrow (V \otimes W)^{h\Gamma \rel}
\end{aligned}
\end{equation}
which induce \eqref{eq:relative-product}. These maps inherit the (strict) associativity property. 
We have one more condition, which concerns commutativity of \eqref{eq:chain-exterior-product}. This does not have to hold on the nose, but we want to have homotopy commutativity, and the underlying chain homotopy should be compatible with the forgetful maps in \eqref{eq:chain-include-and-forget}. This means that its restriction can be used to ensure the homotopy commutativity of \eqref{eq:chain-relative-product}.

Explicitly, we adopt the most classical chain level model, which has all the properties discussed above. It is given by the bar construction
\begin{equation} \label{eq:bar-construction}
V^{h\Gamma} \stackrel{\mathrm{def}}{=} \prod_{k \geq 0} \mathit{Map}(\Gamma^k,V)[-k].
\end{equation}
One can think of it as using the bar resolution of $\bF$ to compute \eqref{eq:derived-hom}; or equivalently, as the space of maps $\bF \rightarrow V$ in the dg category of $A_\infty$-modules over $\bF\Gamma$. Explicitly, the resulting differential is
\begin{equation} \label{eq:bar-differential}
\begin{aligned}
& (d \phi)^k(g_1,\dots,g_k) = (-1)^k d (\phi^k(g_1,\dots,g_k)) + 
g_1 \phi^{k-1}(g_2,\dots,g_k) \\ & \qquad + \sum_j (-1)^j \phi^{k-1}(g_1,\dots,g_jg_{j+1},\dots,g_k) + (-1)^k\phi^{k-1}(g_1,\dots,g_{k-1}).
\end{aligned}
\end{equation}
This comes with obvious maps \eqref{eq:chain-include-and-forget}. The subcomplex \eqref{eq:kernel} consists of those $\phi$ with $\phi^0 = 0$. The product \eqref{eq:chain-exterior-product} is
\begin{equation} \label{eq:chain-exterior-product-2}
\begin{aligned}
& (\phi \times \psi)^k(g_1,\dots,g_k) 
\\ & = \sum_j (-1)^{(|\phi|-j)(k-j)} \phi^j(g_1,\dots,g_j) \otimes
g_1\dots g_j \cdot \psi^{k-j}(g_{j+1},\dots,g_k).
\end{aligned}
\end{equation}
For instance,
\begin{equation} \label{eq:example-product}
\begin{aligned}
& (\phi \times \psi)^0 = \phi^0 \otimes \psi^0, \\
& (\phi \times \psi)^1(g) = \phi^1(g) \otimes g\psi^0 + (-1)^{|\phi|} \phi^0 \otimes \psi^1(g).
\end{aligned}
\end{equation}
How does \eqref{eq:chain-exterior-product-2} arise? Let's write $B\Gamma$ for the bar resolution of $\bF$. Bar resolutions have a multiplicative structure, consisting in this case of a canonical map of $\bF \Gamma^2$-modules
\begin{equation} \label{eq:bar-map-1}
B(\Gamma^2) \longrightarrow B\Gamma \otimes B\Gamma,
\end{equation}
see e.g.\ \cite[Equation X.7.2]{maclane}. They are also naturally functorial, which one applies to get a map linear over the diagonal $\bF \Gamma \rightarrow \bF \Gamma^2$,
\begin{equation} \label{eq:bar-map-2}
B\Gamma \longrightarrow B(\Gamma^2).
\end{equation}
One applies $\mathit{hom}_{\bF \Gamma^2}(\cdot, V \otimes W)$ to \eqref{eq:bar-map-1}, and similarly with \eqref{eq:bar-map-2}, to obtain \eqref{eq:chain-exterior-product-2}. This in particular shows that on the level of cohomology, it reproduces the derived category construction \eqref{eq:exterior-product-2}. The homotopy commutativity of \eqref{eq:chain-exterior-product-2} can be proved by a version of the classical Steenrod $\smile_1$-product. That chain homotopy is compatible with forgetful maps, as desired. We omit the details.

\subsection{Topological groups\label{subsec:topological-groups}}
For a topological group $G$, the chains $C_{-*}(G)$ form a (unital) differential graded algebra, which comes with an augmentation
\begin{equation} \label{eq:augmentation-2}
\epsilon: C_{-*}(G) \longrightarrow \bF.
\end{equation}
By a chain complex $V$ with a $G$-action, we mean a differential graded module over $C_{-*}(G)$. The simplest example is again the trivial module $V = \bF$, coming from \eqref{eq:augmentation-2}. If $V$ and $W$ carry $G$-actions, so does their tensor product (over $\bF$) $V \otimes W$. Conceptually, this is because $V \otimes W$ is naturally a $C_{-*}(G)^{\otimes 2}$-module, and from that one gets a $C_{-*}(G)$-module structure through the diagonal map. Diagrammatically, this means that the $G$-action on $V \otimes W$ is defined by
\begin{equation} \label{eq:tensor-product}
\begin{aligned}
& C_{-*}(G) \otimes (V \otimes W) \xrightarrow{\text{diagonal}} (C_{-*}(G) \otimes C_{-*}(G)) \otimes (V \otimes W) 
\\ & \quad \iso (C_{-*}(G) \otimes V) \otimes (C_{-*}(G) \otimes W) \xrightarrow{\text{$G$-actions}} V \otimes W,
\end{aligned}
\end{equation}
where the diagonal is the standard Alexander-Whitney map. Explicitly, 
\begin{equation} \label{eq:tensor-product-2}
c \cdot (v \otimes w) = \sum_i (-1)^{(|c|-i)|v|} (c^{\leq i} \cdot v) \otimes (c^{\geq i} \cdot w).
\end{equation}
Here $c$ is a simplex in $C_{-*}(G)$ of degree $|c|$; and $c^{\leq i}$, $c^{\geq i}$ are the faces of $c$ which occur in the Alexander-Whitney diagonal. The tensor product defined in this way is strictly associative, and is commutative up to quasi-isomorphism. If $\bF$ is given the trivial action, we have $V \otimes \bF = \bF \otimes V = V$.

One defines equivariant cohomology and its relative version as before,
\begin{equation} \label{eq:derived-hom-2}
\begin{aligned}
& H^*_G(V) = \mathit{DHom}_{C_{-*}(G)}(\bF,V), \\
& H^*_{G\rel}(V) = \mathit{DHom}_{C_{-*}(G)}(\mathit{ker}(\epsilon)[1],V).
\end{aligned}
\end{equation}
These groups come with maps \eqref{eq:include-and-forget} (one defines $V^G$ as the subcomplex on which each element of $G$ acts as the identity, and all higher simplices act as zero), as well as with exterior products. Finally, for manifolds, \eqref{eq:poincare-duality} continues to hold.

\begin{remark} \label{th:shortcoming-0}
Take topological spaces $X$ and $Y$ on which $G$ acts. Then, $C_{-*}(X \times Y)$ has an obvious $G$-action, and so does $C_{-*}(X) \otimes C_{-*}(Y)$ via \eqref{eq:tensor-product-2}.
However, \eqref{eq:cross-product} is compatible with those actions only in a homotopical sense. That creates a problem for writing down a chain map underlying \eqref{eq:topological-exterior-product}. In our applications, we will bypass this by replacing the topological group by a discrete one, which gives rise to the same equivariant cohomology groups (Section \ref{subsec:p-circle}).
A closely related issue arises even before we get to topology: the tensor product $G$-action \eqref{eq:tensor-product-2} is commutative only up to quasi-isomorphism. As a consequence, commutativity of the exterior product on equivariant cohomology is a little harder to establish. For the instances relevant to us, one can again evade that extra work by using a discrete group instead.
\end{remark}

On the chain level, the bar construction is now
\begin{equation} \label{eq:bar-construction-2}
V^{hG} = \prod_{k \geq 0} \mathit{Hom}(C_{-*}(G)^{\otimes k}, V)[-k].
\end{equation}
with differential
\begin{equation} \label{eq:bar-differential-2}
\begin{aligned}
& (d\phi)^k(c_1,\dots,c_k) \\ & \quad = (-1)^k d(\phi^k(c_1,\dots,c_k)) 
- \sum_i (-1)^{|\phi|+|c_1|+\cdots+|c_{i-1}|} \phi^k(c_1,\dots, dc_i,\dots,c_k) 
\\ & \quad 
+ (-1)^{(|\phi|-k+1)|c_1|} c_1 \cdot \phi^{k-1}(c_2,\dots,c_k) 
+ \sum_i (-1)^i \phi^{k-1}(c_1,\dots,c_ic_{i+1},\dots,c_k) 
\\ & \quad
+ (-1)^k \phi^{k-1}(c_1,\dots,c_{k-1})\epsilon(c_k).
\end{aligned}
\end{equation}
Here, $|\phi|$ is the degree in \eqref{eq:bar-construction-2}. For instance,
\begin{equation} \label{eq:example-differential}
\begin{aligned}
& (d\phi)^0 = d(\phi^0), \\
& (d\phi)^1(c) = -d(\phi^1(c)) - (-1)^{|\phi|} \phi^1(dc)
+(-1)^{|\phi|\,|c|} c \cdot\phi^0 - \phi^0 \,\epsilon(c).
\end{aligned}
\end{equation}
The appropriate generalization of \eqref{eq:chain-exterior-product-2} is
\begin{equation} \label{eq:chain-exterior-product-3}
\begin{aligned}
& (\phi \times \psi)^k(c_1,\dots,c_k) = 
\\ & \sum_{j, i_1,\dots,i_j}(-1)^\S  \phi^j(c_1^{\leq i_1},\dots,c_j^{\leq i_j}) \otimes
(c_1^{\geq i_1} \cdots c_j^{\geq i_j}) \cdot \psi^{k-j}(c_{j+1},\dots,c_k)
\end{aligned}
\end{equation}
for simplices $c_1,\dots,c_k$ (in the term after the tensor product, we take the product of certain faces of those simplices in the ring $C_{-*}(G)$, and then apply that to $\psi^{k-j}(\cdots) \in W$ using the module structure). Finally, the sign is the parity of
\begin{equation}
\S = (|\phi|-j)(k-j) + (|\psi|+k-j) (|c_1|+\cdots+|c_j|) + \sum_{p<q\leq j} |c_p^{\geq i_p}|\,|c_q^{\leq i_q}|.
\end{equation}
For instance,
\begin{equation} \label{eq:example-product-2}
\begin{aligned}
& (\phi \times \psi)^0 = \phi^0 \otimes \psi^0, \\
& (\phi \times \psi)^1(c) = \sum_i (-1)^{|\psi|\,|c|} \phi^1(c^{\leq i}) \otimes c^{\geq i} \psi^0 \; + (-1)^{|\phi|}
\phi^0 \otimes \psi^1(c).
\end{aligned}
\end{equation}

\section{Specific groups}
As mentioned in the introduction, our application of equivariant cohomology involves only two particularly simple cases: the circle (as a topological group), and its prime order finite subgroups. 
We introduce a third group as a discrete replacement for the circle (this is a well-known idea, see e.g.\ \cite{nikolaus-scholze18}).

\subsection{The circle}
Take the topological group $G = S^1 = \bR/\bZ$. This has
\begin{equation} \label{eq:s1-cohomology}
H^*_G(\bF) = \bF[t], \;\; |t| = 2.
\end{equation}
In terms of \eqref{eq:bar-construction-2}, the generator is represented by a cocycle in $\bF^{hG}$ with two components:
\begin{equation} \label{eq:t-cocycle}
\begin{aligned}
& t^1: C_1(G) \rightarrow \bF, && t^1(c) = \lfloor \tilde{c}(0) \rfloor - \lfloor \tilde{c}(1) \rfloor, \\
& t^2: C_0(G)^{\otimes 2} \rightarrow \bF, 
&& t^2(g_1,g_2) = \lfloor \tilde{g}_1+ \tilde{g}_2 \rfloor - \lfloor \tilde{g}_1 \rfloor - \lfloor \tilde{g}_2 \rfloor.
\end{aligned}
\end{equation}
Here, $\tilde{c}: [0,1] \rightarrow \bR$ is a lift of $c: [0,1] \rightarrow G$; similarly, $g_i \in G$, and the $\tilde{g}_i$ are preimages in $\bR$. One way to think of \eqref{eq:t-cocycle} is by starting with the degree $1$ cochain $\tilde{\phi}^1(\tilde{g}) = -\lfloor \tilde{g} \rfloor$ for the group $\tilde{G} = \bR$; then $\tilde{t} = d\tilde{\phi}$ descends to our cocycle $t$ on $G$, but is no longer exact there.

Let's look at the relative version of the theory. The coefficients are obviously
\begin{equation} \label{eq:relative-g}
H^*_{G \rel}(\bF) = t\bF[t] \subset H^*_G(\bF).
\end{equation}
In particular, multiplication by $t$ can be viewed as an instance of \eqref{eq:relative-product}, which means that it is naturally a map
\begin{equation} \label{eq:multiplication-by-t}
t: H^*_G(V) \longrightarrow H^{*+2}_{G\rel}(V).
\end{equation}

\begin{lemma} \label{th:relative-t}
(i) The map \eqref{eq:multiplication-by-t} is an isomorphism.

(ii) Consider the composition
\begin{equation} \label{eq:three-maps}
H^*(V) \longrightarrow H^{*+1}_{G \rel}(V) \stackrel{\iso}{\longleftarrow} H^{*-1}_G(V) \longrightarrow H^{*-1}(V),
\end{equation}
of: the connecting map from \eqref{eq:relative-les}; the inverse of \eqref{eq:multiplication-by-t}; and the forgetful map. This composition equals loop rotation, meaning the action of $[G]$ on the $H_{-*}(G)$-module $H^*(V)$.
\end{lemma}

\begin{proof}
(i) The map \eqref{eq:multiplication-by-t} is functorial in $V$. In view of \eqref{eq:derived-hom-2} and the Yoneda Lemma, it must be given by some degree $2$ element of
\begin{equation} \label{eq:map-from-relative}
D\mathit{hom}_{C_{-*}(G)}(\mathit{ker}(\epsilon)[1],\bF).
\end{equation}
The cohomology of $\mathit{ker}(\epsilon)$ is one-dimensional, concentrated in degree $-1$. One can see (either by direct construction, or by an obstruction theory argument) that $\mathit{ker}(\epsilon)$ is quasi-isomorphic to $\bF[1]$ as a $C_{-*}(G)$-module, so that \eqref{eq:map-from-relative} is isomorphic to $H^{*-2}_G(\bF)$. Any degree $2$ element of that space is either zero, or an isomorphism (in the derived category). Since \eqref{eq:multiplication-by-t} is nonzero for $V = \bF$, only the second of those options is possible. 

(ii) 
Consider the composition of the first two maps in \eqref{eq:three-maps}. This is functorial in $V$, hence must be given by some degree $-1$ element of
\begin{equation}
D\mathit{hom}_{C_{-*}(G)}(\bF,C_{-*}(G)) = H^*_G(C_{-*}(G)).
\end{equation}
By \eqref{eq:poincare-duality} this group is one-dimensional, and its generator maps to $[G]$ under the forgetful map. Hence, the desired equality is true up to a universal multiplicative constant. 

To determine that constant, it is sufficient to consider the case $V = C_{-*}(G)$. If $\phi \in V^{hG \rel}$ is a cocycle of degree $1$, then $\phi^1: C_*(G) \rightarrow C_*(G)$ is a chain map, see \eqref{eq:example-differential}. The induced map on $H_1(G)$ gives a well-defined map (the topological degree; this has nothing to do with gradings)
\begin{equation} \label{eq:degree}
\mathrm{Deg}: H^1(V^{hG \rel}) \longrightarrow \bF.
\end{equation}
Start with the neutral element $e \in G \subset C_0(G)$. Again by \eqref{eq:example-differential}, its image under the connecting map $H^0(V) \rightarrow H^1_{G\rel}(V)$ is represented by a cocycle $\phi$ with
\begin{equation}
\phi^1(c) = c - e \, \epsilon(c),
\end{equation}
hence $\mathrm{Deg}(\phi) = 1$. On the other hand, take any degree $-1$ cocycle $\psi \in V^{hG}$ whose constant term satisfies $[\psi^0] = [G] \in H_1(G)$. If we take its product with the cocycle representative for $t$ from \eqref{eq:t-cocycle}, then by \eqref{eq:example-product-2},
\begin{equation}
(t \times \psi)^1(c) =- t^1(c) \, (c^{\geq 1} \psi^0) \;\; \text{ for } c \in V^{-1} = C_1(G),
\end{equation}
which maps to $1$ under \eqref{eq:degree}. Hence, the image of $[e]$ under the connecting homomorphism agrees with $t[\psi]$, as desired.
\end{proof}

\subsection{Finite cyclic groups\label{subsec:finite-cyclic-group}}
Let $\Gamma \subset S^1$ be the cyclic group of odd prime order $p$, with generator $\gamma = 1/p$.  We will work with coefficients in $\bF = \bF_p$. One has
\begin{equation} \label{eq:gamma-coefficients}
H^*_\Gamma(\bF) \iso \bF[t,\theta],
\end{equation}
where $|t| =2$ as before, and $|\theta| = 1$; the ring structure is supercommutative, meaning that $t\theta = \theta t$ and $\theta^2 = 0$. The relative version is correspondingly
\begin{equation}
H^*_{\Gamma\rel}(\bF) = \theta \bF[t] \oplus t\bF[t].
\end{equation}
On the complex $\bF^{h\Gamma}$, the class $\theta$ is represented by the obvious map
\begin{equation}
\label{eq:theta-cocycle}
\theta^1: \Gamma \rightarrow \bF, \;\; \theta^1(\gamma^i) = i.
\end{equation}
The class $t$ is represented by the restriction of \eqref{eq:t-cocycle}, which can be explicitly written as
\begin{equation}
\label{eq:t-cocycle-2}
t^2: \Gamma^2 \rightarrow \bF, \;\;
t^2(\gamma^i,\gamma^j) = \lfloor (i+j)/p \rfloor \text{ for $0 \leq i,j < p$.}
\end{equation}
Consider the short exact sequence
\begin{equation} \label{eq:witt-vectors}
0 \rightarrow \bF \longrightarrow \tilde{\bF} = (\bZ/p^2) \longrightarrow \bF \rightarrow 0,
\end{equation}
From the analogue of $\bF^{h\Gamma}$ with $\tilde{\bF}$-coefficients, one gets the Bockstein operation $\beta$ on $H^*_{\Gamma}(\bF)$. By looking at \eqref{eq:theta-cocycle} and \eqref{eq:t-cocycle-2}, one sees that
\begin{equation} \label{eq:beta-theta}
\beta(\theta) = t.
\end{equation}

\begin{remark} \label{th:ainfinity}
There is another relation between $\theta$ and $t$, via chain level multiplicative structures, which can be expressed as follows: a minimal $A_\infty$-model for $\bF^{h\Gamma}$ is given by \eqref{eq:gamma-coefficients} with the nontrivial higher order $A_\infty$-operation
\begin{equation} \label{eq:mu-p}
\mu^p(\theta,\dots,\theta) = \pm t
\end{equation}
One can generalize this model to the $\Gamma$-equivariant cohomology of any complex $V$. The generalization is $V[[t]][\theta]$, with the well-known $t$-linear differential
\begin{equation}
\begin{aligned}
& v \mapsto  dv + \theta (\gamma v - v), \\
& \theta \bar{v}  \mapsto -\theta\, d\bar{v} + t(\bar{v} + \gamma \bar{v} + \cdots + \gamma^{p-1}\bar{v})
\end{aligned}
\end{equation}
($d$ being the differential on $V$). This time, there are nontrivial ``exterior $A_\infty$-products''
%
%
\begin{equation}
\mu^n: V_1[[t]][\theta] \otimes \cdots \otimes V_n[[t]][\theta] \longrightarrow (V_1 \otimes \cdots \otimes V_n)[[t]][\theta]
\end{equation}
of all orders $2 \leq n \leq p$, which are again $t$-linear and satisfy
\begin{equation}
\mu^n(\theta v_1,\dots,\theta v_n) = \pm t \!\!\!\! \sum_{0 \leq i_1 < \cdots < i_n < p}
\!\! \gamma^{i_1}v_1 \otimes \cdots \otimes \gamma^{i_n} v_n;
\end{equation}
in the special case $V_1 = \cdots = V_n = \bF$, all of them cancel except \eqref{eq:mu-p}.
Since we're not using this model, we will not explain how to derive it (or clarify the signs, which depend on the conventions used for $A_\infty$-structures). Morally, the upshot is that the multiplicative structure of the $\Gamma$-equivariant theory is a little more complicated than it might appear from \eqref{eq:gamma-coefficients}.
\end{remark}

Take a chain complex $V$ with a $G$-action. The restriction map in equivariant cohomology,
\begin{equation} \label{eq:g-to-gamma}
H^*_G(V) \longrightarrow H^*_{\Gamma}(V),
\end{equation}
is implemented in our chain level context by restricting cochains in \eqref{eq:bar-construction-2} from $C_{-*}(G)$ to points in $\Gamma$. For $V = \bF$, a comparison of \eqref{eq:t-cocycle} and \eqref{eq:t-cocycle-2} shows that this maps $t \in H^2_G(\bF)$ to the corresponding element in $H^2_\Gamma(\bF)$, so that it yields an isomorphism $H^*_G(\bF)[\theta] \iso H^*_\Gamma(\bF)$. More generally:

\begin{lemma} \label{th:comparison}
Let $V$ be a complex with $G$-action. The restriction map \eqref{eq:g-to-gamma} is injective, and induces an isomorphism of $H^*_\Gamma(\bF)$-modules,
\begin{equation} \label{eq:theta-isomorphism}
H^*_G(V)[\theta] \stackrel{\iso}{\longrightarrow} H^*_{\Gamma}(V).
\end{equation}
\end{lemma}

\begin{proof}
One uses adjunction:
\begin{equation} \label{eq:adjunction}
\begin{aligned} &
H^*_{\Gamma}(V) \iso D\mathit{Hom}_{\bF\Gamma}(\bF, V) \\ & \qquad \iso
D\mathit{Hom}_{C_{-*}(G)}(C_{-*}(G) \otimes_{\bF \Gamma} \bF, V) 
= D\mathit{Hom}_{C_{-*}(G)}(C_{-*}(G/\Gamma), V).
\end{aligned}
\end{equation}
On general grounds, the $\otimes$ in the formula should be a derived one; but $C_{-*}(G)$ being a bounded above complex of free $\bF\Gamma$-modules, that is unnecessary. 

The augmentation
\begin{equation} \label{eq:adjunction-2}
C_{-*}(G) \otimes_{\bF \Gamma} \bF \xrightarrow{\epsilon \otimes 1} \bF.
\end{equation}
induces a map $D\mathit{Hom}_{C_{-*}(G)}(\bF,V) \rightarrow D\mathit{Hom}_{C_{-*}(G)}(
C_{-*}(G) \otimes_{\bF \Gamma} \bF,V)$, which under the identification \eqref{eq:adjunction} is the restriction \eqref{eq:g-to-gamma}. Similarly, we have a map (in the derived category)
\begin{equation} \label{eq:adjunction-3}
C_{-*}(G) \otimes_{\bF\Gamma} \bF \xrightarrow{1 \otimes [\theta]} C_{-*}(G)[1] \otimes_{\bF\Gamma} \bF \xrightarrow{\eqref{eq:adjunction-2}} \bF[1],
\end{equation}
which induces the $\theta$-component of \eqref{eq:theta-isomorphism}. On cohomology, the first map in \eqref{eq:adjunction-3} is the action of $\theta$ on $H_{-*}^\Gamma(G)$. Equivalently, the $\Gamma$-covering determines a class in $H^1(G/\Gamma)$, and we're looking at its cap action on $H_{-*}(G/\Gamma)$. Clearly, this action is nontrivial. The consequence is that the combination of \eqref{eq:adjunction-2} and \eqref{eq:adjunction-3} gives a quasi-isomorphism $C_{-*}(G) \otimes_{\bF\Gamma} \bF \rightarrow \bF \oplus \bF[1]$. This, by adjunction, implies that \eqref{eq:theta-isomorphism} is an isomorphism.
\end{proof}

\subsection{The $p$-circle\label{subsec:p-circle}}
We continue with $p$ and $\bF$ as before. Write $\Gamma_N \subset G$ for the subgroup of elements of order dividing $p^N$ (so that $\Gamma_1 = \Gamma$). Let $\Pi$ be the union of all those subgroups, still considered as a discrete group. We call this the $p$-circle, but in the literature it is usually called the Pr\"ufer group, and written as $\bZ/p^\infty$. It satisfies
\begin{equation} \label{eq:p-circle-cohomology}
H^*_\Pi(\bF) = \bF[t].
\end{equation}
The generator $t$ is the image of its counterpart for $G$ under the restriction map. On the chain level, it retains the component $t^2: \Pi^2 \rightarrow \bF$ of \eqref{eq:t-cocycle}. One can view this as a special case of the following:

\begin{lemma} \label{th:comparison-2}
For any complex $V$ with a $G$-action, 
\begin{equation} \label{eq:discretize}
H^*_G(V) \stackrel{\iso}{\longrightarrow} H^*_{\Pi}(V).
\end{equation}
\end{lemma}

\begin{proof}
This uses adjunction, as in Lemma \ref{th:comparison}. The main step is to look at the $C_{-*}(G)$-module
\begin{equation}
C_{-*}(G) \otimes_{\bF\Pi} \bF = \underrightarrow{\lim}_N \, C_{-*}(G) \otimes_{\bF\Gamma_N} \bF = \underrightarrow{\lim}_N \, C_{-*}(G/\Gamma_N).
\end{equation}
The augmentation map from this to $\bF$ is easily seen to be a quasi-isomorphism.
\end{proof}

The combination of Lemmas \ref{th:comparison} and \ref{th:comparison-2} shows that for any complex $V$ with a $G$-action,
\begin{equation} \label{eq:discretize-2}
H^*_{\Pi}(V)[\theta] \stackrel{\iso}{\longrightarrow} H^*_{\Gamma}(V).
\end{equation}

\begin{remark} \label{th:prufer}
The isomorphism \eqref{eq:discretize-2} fails for general $\Pi$-actions, see e.g.\ $V = \bF\Pi$.
\end{remark}

Occasionally, we will encounter arguments involving the larger coefficient ring $\tilde{\bF}$. The context is that $V$ can be lifted to a complex $\tilde{V}$ of free $\tilde{\bF}$-modules, with the corresponding notion of $G$-action. The proof of Lemma \ref{th:comparison-2} goes through in that context, so that the well-known behaviour of $G$-equivariant cohomology with $\tilde{\bF}$-coefficients also carries over to $\Pi$. Returning to $\bF$-coefficients, the choice of such a lift $\tilde{V}$ gives rise to Bockstein operators $\beta$, which sit in a commutative diagram
\begin{equation}
\xymatrix{
H^*_G(V) \ar[d]_-{\iso} \ar[rr]^-{\beta} && H^{*+1}_G(V) \ar[d]^-{\iso} \\
H^*_\Pi(V) \ar[rr]^-{\beta} && H^{*+1}_\Pi(V).
}
\end{equation}

\section{The TQFT framework\label{sec:tqft}}
In our framework, operations are parametrized by cylinders with marked points. More specifically, the parameter spaces are symmetric products of the cylinder. We apply $G$-equivariant cohomology to the complex of singular chains on these spaces, where $G = S^1$ rotates the cylinder; on the chain level, that is implemented using $\Pi$ instead of $G$. Having set up the general theory, we use it to first define a $q$-deformed differential, and then abstract analogues of Steenrod operations. 

\subsection{Parameter spaces\label{subsec:parameter-spaces}}
Consider the cylinder $\bR \times S^1$, with coordinates $z = (s,t)$. For any integer $m > 0$, define
\begin{equation} \label{eq:p-space}
\frakP_m = \mathit{Sym}_m(\bR \times S^1)/\bR,
\end{equation}
where $\bR$ acts by $s$-translations. We usually write elements of $\mathit{Sym}_m(\bR \times S^1)$ as $\Sigma =  (z_1,\dots,z_m)$, where $z_k = (s_k,t_k)$, and their equivalence classes with respect to translation as $[\Sigma] = [z_1,\dots,z_m]$; however, one should keep in mind that these are really unordered collections. The spaces \eqref{eq:p-space} come with natural compactifications, which allow points of $\Sigma$ to separate out to $\pm\infty$, and the cylinder accordingly to split into several ones. The compactifications are manifolds with corners $\bar\frakP_m$, whose strata are labeled by partitions:
\begin{equation}
\bar\frakP_m = \!\! \coprod_{\substack{l>0, \, m^1,\dots,m^l > 0 \\ m^1+\cdots+m^l = m}} \!\!
\frakP_{m^1} \times \cdots \times\frakP_{m^l}.
\end{equation}
In the simplest case, $\frakP_1 = \bar\frakP_1 = S^1$. Stacking cylinders next to each other defines associative composition maps
\begin{equation} \label{eq:gluing}
\bar\frakP_m \times \bar\frakP_n \longrightarrow \bar\frakP_{m+n}.
\end{equation}
It is convenient to add another space
\begin{equation}
\frakP_0 = \bar\frakP_0 = \mathit{point},
\end{equation}
extending \eqref{eq:gluing} so that left and right composition with $\bar\frakP_0$ is the identity.
The spaces $\frakP_m$ and their compactifications carry actions of $G$, rotating the $t$-factor of the cylinder, which are compatible with \eqref{eq:gluing}. 

So far, the discussion would apply equally to other versions of configuration spaces on the cylinder. We now analyze the topology, where symmetric products turn out to be particularly simple. Namely, for each $m>0$, the map 
\begin{equation} \label{eq:morton}
\begin{aligned}
& \frakP_m \stackrel{\htp}{\longrightarrow} \mathit{Sym}_m(S^1) \stackrel{\htp}{\longrightarrow} S^1, \\
& [(s_1,t_1),\dots, (s_m,t_m)] \longmapsto
t_1 + \cdots + t_m,
\end{aligned}
\end{equation}
and its continuous extension to $\bar\frakP_m$, are homotopy equivalences \cite{morton67}. If one wants to make this map $G$-equivariant, the target should be equipped with the action which rotates $G$ with multiplicity $m$. 

Throughout, we work over $\bF = \bF_p$ for odd $p$. Consider
\begin{equation} \label{eq:p}
P = \bigoplus_{m \geq 0} P_m, \quad P_m = C_{-*}(\bar\frakP_m).
\end{equation}
We have reversed the sign of the grading on chains, so as to make it agree with general cohomological conventions for algebraic structures. The composition maps \eqref{eq:gluing} induce chain maps
\begin{equation}
P_i \otimes P_j \longrightarrow P_{i+j},
\end{equation}
which turn $P$ into a bigraded dga. The unit element is $e = \mathit{point} \in P_0^0 = C_0(\bar\frakP_0)$. On the level of cohomology, \eqref{eq:morton} implies the following:

\begin{lemma} \label{th:cohomology}
As a bigraded algebra,
\begin{equation} \label{eq:iota-lambda}
\begin{aligned}
& H^*(P) \iso \bF[\iota,\lambda], \\
& \iota = [\mathit{point}] \in H^0(P_1) = H_0(\frakP_1), \\
& \lambda = [S^1] \in H^{-1}(P_1) = H_1(\frakP_1).
\end{aligned}
\end{equation}
In the notation, supercommutativity is implicit, in particular $\lambda^2 = 0$.
\end{lemma}

We now consider the equivariant flavour of the theory. The most obvious version is
\begin{equation} \label{eq:pg}
P^{hG} = \bigoplus_{m \geq 0} P_m^{hG},
\end{equation}
which has 
\begin{equation} \label{eq:pg-cohomology}
H^*(P^{hG}) = \bigoplus_{m \geq 0} H^*(P_m^{hG}) = \bigoplus_{m \geq 0} H^*_G(P_m).
\end{equation}
However, there are issues in defining a multiplicative structure on \eqref{eq:pg}, see Remark \ref{th:shortcoming-0}. We can use the discrete replacement instead,
\begin{equation} \label{eq:ppi}
P^{h\Pi} = \bigoplus_{m \geq 0} P_m^{h\Pi},
\end{equation}
which has a (bigraded) dga structure, defined by combining \eqref{eq:gluing} with \eqref{eq:chain-exterior-product}. On cohomology, there is no difference between \eqref{eq:pg} and \eqref{eq:ppi}, thanks to Lemma \ref{th:comparison-2}. In practice, this means that for cohomology level arguments \eqref{eq:pg-cohomology} will be used, since that is more intuitive; while all cochain level constructions requiring multiplicative structures will be based on \eqref{eq:ppi}.

\begin{lemma} \label{th:equivariant-cohomology}
(i) $H^*_G(P_0) \iso H^*_G(\bF) = \bF[t]$.

(ii) If $m$ is not a multiple of $p$, the map $H^*_G(P_m) \rightarrow H^*(P_m)$ is an isomorphism in degree $* = -1$. In all other degrees, $H^*_G(P_m)$ is zero.

(iii) If $m>0$ is a multiple of $p$, then $H^*_G(P_m)$ is a free graded module over $\bF[t]$, with one generator each in degrees $-1$ and $0$. The map $H^*_G(P_m) \rightarrow H^*(P_m)$ is an isomorphism in degrees $* = -1,0$, so that one gets an isomorphism 
\begin{equation} \label{eq:p-splitting}
H^*_G(P_m) \iso H^*(P_m)[t].
\end{equation}

(iv) For any $m$, $H^*_{G \rel}(P_m)$ is concentrated in degrees $\geq 1$; and if $m$ is not a multiple of $\mathrm{char}(\bF)$, then it is concentrated in degree $1$.
\end{lemma}

\begin{proof}
(i) is obvious. For (ii) and (iii), we can use \eqref{eq:morton} and \eqref{eq:poincare-duality}. This gives isomorphisms
\begin{equation} \label{eq:lens-space}
H^*_G(P_m) \iso H^{*+1}_G(S^1) \iso H^{*+1}(L_m).
\end{equation}
Here, $G$ acts on $S^1$ via $m$-fold rotation; and the infinite-dimensional lens space $L_m = S^\infty/(\bZ/m)$ is the result of applying the Borel construction. The action of $t$ on the right hand side comes from the structure of $H^*(L_m)$ as a module over $H^*(\bC P^\infty)$. From there, the desired computation is straightforward. Finally, (iv) follows from the previous parts by using the long exact sequence \eqref{eq:relative-les}.
\end{proof}

For any $m>0$, write 
\begin{equation} \label{eq:equivariant-generator}
 (\lambda\iota^{m-1})_G \in H^{-1}_G(P_m)
\end{equation}
for the unique lift of $\lambda\iota^{m-1}$.  In case (iii) of Lemma \ref{th:equivariant-cohomology}, we similarly write
\begin{equation} \label{eq:equivariant-generator-2}
\iota^m_G \in H^0_G(P_m)
\end{equation}
for the unique lift of $\iota^m$. 

We will need an additional piece of information. $P_m^{hG}$ can be defined over $\tilde{\bF}$, giving rise to a Bockstein for the $\bF$-coefficient version,
\begin{equation} \label{eq:bockstein}
\beta: H^*_G(P_m) \longrightarrow H^{*+1}_G(P_m).
\end{equation}

\begin{lemma} \label{th:bockstein}
Suppose that $m>0$ is a multiple of $p$. Then the classes \eqref{eq:equivariant-generator}, \eqref{eq:equivariant-generator-2} satisfy
\begin{equation}
\beta(\iota^m_G) = (m/p)\,t (\lambda\iota^{m-1})_G \in H^1_G(P_m).
\end{equation}
\end{lemma}

\begin{proof}
We can replace $P_m$ with the homotopy equivalent $V = C_{-*}(S^1)$, carrying the $m$-fold rotational action. One has a topological degree map as in the (similar, but not identical) context of Lemma \ref{th:relative-t}(ii); and the class $t(\lambda\iota^{m-1})_G$ has topological degree $1$. By definition, one can represent $\iota^m_G$ by a cocycle
\begin{equation}
\phi \in V^{hG}, \;\; \phi^0 = e.
\end{equation}
Choose a lift $\tilde{\phi}$ to a cochain with coefficients in $\tilde{\bF}$, with the same constant term. By \eqref{eq:example-differential}, 
\begin{equation}
\begin{aligned}
& (d\tilde{\phi})^0 = 0, \\
& (d\tilde{\phi})^1([G]) = -d(\tilde{\phi}^1([G])) + [G] \cdot e
\end{aligned}
\end{equation}
for some fundamental chain $[G]$. Because of the specific action of $G$ on $S^1$ we are considering, $[G] \cdot e$ represents $m$ times the fundamental class of the circle. After dividing by $p$, one sees that $\beta(\iota^m_G)$ is represented by an element of $V^{hG \rel}$ which has degree $m/p$.
\end{proof}

\subsection{The deformed complex}
We now introduce a formal parameter $q$.

\begin{lemma} \label{th:deltaq}
(i) There is a series
\begin{equation} \label{eq:deltaq}
\begin{aligned}
& \delta_q = \delta_1 q + \delta_2 q^2 + \cdots, \;\; \delta_m \in P_m^{h\Pi,1-2m}, \\
& d\delta_q + \delta^2_q = 0, \\
& [\delta_1] = \lambda_G \in H^{-1}_G(P_1),\;\; \text{ see \eqref{eq:equivariant-generator} for notation.}
\end{aligned}
\end{equation}

(ii) Given $\delta_q$ and $\delta_q'$ both having the properties above, there is
\begin{equation} \label{eq:phiq}
\begin{aligned}
& \phi_q = \phi_1 q+ \phi_2 q^2 + \cdots, \;\; \phi_m \in P_m^{h\Pi,-2m}, \\
& d\phi_q + \delta_q' (e+\phi_q) - (e+\phi_q) \delta_q = 0. 
\end{aligned}
\end{equation}
\end{lemma}

\begin{proof}
(i) We proceed order-by-order in $q$. Having chosen $\delta_1,\dots,\delta_{m-1}$, the obstruction to finding $\delta_m$ lies in $H_G^{2-2m}(P_m)$; but those groups are zero for all $m \geq 2$, by Lemma \ref{th:equivariant-cohomology}.

(ii) The first order obstruction is $[\delta_1' - \delta_1] \in H_G^{-1}(P_1)$, but that is zero by assumption on those leading terms. The remaining obstructions lie in 
$H_G^{1-2m}(P_m)$ for $m \geq 2$, and those groups are zero.
\end{proof}

Take $i,j \in \bZ$. Consider series 
\begin{equation} \label{eq:xq}
x = \sum_{m \geq i} x_m q^{m-i} = x_i + x_{i+1} q + x_{i+2} q^2 + \cdots, \quad x_m \in P_m^{h\Pi,j+2(i-m)}.
\end{equation}
The space of such $x$ is denoted by $Q^{h\Pi,j}_i$. Note that
\begin{equation} \label{eq:multiply-by-q}
qQ^{h\Pi,j}_i \subset Q^{h\Pi,j+2}_{i-1}.
\end{equation}
For instance, $\delta_q \in Q^{h\Pi,1}_0$ (with the additional property that the $q$-constant term is zero). Having chosen $\delta_q$, we equip $Q^{h\Pi}$ with the product coming from $q$-linearly extending that on $P^{h\Pi}$, as well as with the deformed differential
\begin{equation} \label{eq:dq}
d_q x \stackrel{\mathrm{def}}{=} dx + \delta_q x - (-1)^{|x|} x \delta_q.
\end{equation}
This sends $Q^{h\Pi,*}_i$ to $Q^{h\Pi,*+1}_i$, so that the whole of $Q^{h\Pi}$ is again a bigraded dga. There is also a non-equivariant version $Q$ of this structure, which consists of series with coefficients in $P$. The differential is defined as in \eqref{eq:dq}, using the image of $\delta_q$ under the forgetful maps. This comes with a forgetful map (of bigraded dga's)
\begin{equation} \label{eq:q-forget}
Q^{h\Pi} \longrightarrow Q.
\end{equation}
Finally, one can define a relative version $Q^{h\Pi \rel}$ as in \eqref{eq:kernel}, which is a two-sided ideal in $Q^{h\Pi}$. 

The role of Lemma \ref{th:deltaq}(ii) in this context is to provide a uniqueness theorem (up to quasi-isomorphism) for $Q^{h\Pi}$, following a familiar strategy. One considers the product $Q^{h\Pi} \otimes I$, where $I$ is the cochain dg algebra of the interval:
\begin{equation}
\begin{aligned}
& I = \bF u \oplus \bF u' \oplus \bF f, \;\; |u| = |u'| = 0, \; |f| = 1, \\
& u^2 = u, \;\; (u')^2 = u', \;\; uu' = u'u = 0, \;\; u'f = fu = f, \;\; du' = f = -du.
\end{aligned}
\end{equation}
In the situation of \eqref{eq:phiq}, the element $\delta_q \otimes u + \delta'_q \otimes u' - \phi_q \otimes f$ equips $Q^{h\Pi} \otimes I$ with a deformed dga structure. This comes with two quasi-isomorphic projections to copies of $Q^{h\Pi}$, with the dga structure defined by $\delta_q$ respectively $\delta_q'$ on $Q^{h\Pi}$. For the discussion from now, we will assume that one choice of $\delta_q$ has been fixed.

We will now look at the cohomology of $Q^{h\Pi}$ and its relative counterpart $Q^{h\Pi \rel}$; not in full generality, but selecting specific situations which are important for our construction. This will usually be done by isolating the leading term in \eqref{eq:xq}, which means setting $q = 0$, and showing that the remaining terms do not contribute.

\begin{lemma} \label{th:q-zero-0}
(i) The map (obtained by setting $q=0$) 
\begin{equation} \label{eq:q-zero-0}
H^*(Q^{h\Pi}_0) \longrightarrow H^*_G(P_0) \iso \bF[t]
\end{equation}
is an isomorphism in degrees $\leq 0$, as well as in degree $2$.

(ii) The relative version of that map,
\begin{equation} \label{eq:q-zero-0b}
H^*(Q^{h\Pi \rel}_0) \longrightarrow H^*_{G \rel}(P_0) \iso H^*_{G\rel}(\bF) = t\bF[t],
\end{equation}
is an isomorphism in degrees $\leq 2$.
\end{lemma}

\begin{proof}
Consider the spectral sequence associated to the $q$-adic filtration of $Q^{h\Pi}_0$. This has front page and differential
\begin{equation} \label{eq:q-spectral-sequence}
E_1^{ij} = \begin{cases} 
H^{j-i}_{G}(P_i) & i \geq 0, \\
0 & \text{otherwise}
\end{cases}
\qquad \qquad
d_1(x) = \lambda_G x - (-1)^{|x|} x \lambda_G.
\end{equation}
In total degrees $i+j \leq 3$, the nonzero contributions (from Lemma \ref{th:equivariant-cohomology}) are
\begin{equation} \label{eq:low-degree}
E_1^{0,0} = H^0_G(\bF),  \;\;
E_1^{1,0} = H^{-1}_G(P_1), \;\; 
E_1^{0,2} = H^2_G(\bF), \;\;
E_1^{2,1} = H^{-1}_G(P_2),
\end{equation}
each of which is one-dimensional. In the degrees addressed by (i), all generators lie in the $E_1^{0*}$ column. The generator of $E_1^{0,0}$ is the identity, which clearly survives to $E_\infty$. From that, we immediately get a slightly weaker form of (i): \eqref{eq:q-zero-0} is an isomorphism in degrees $\leq 0$, and injective in degree $2$.

For the relative version of \eqref{eq:q-spectral-sequence}, the analogue of \eqref{eq:low-degree} would be 
\begin{equation} \label{eq:low-degree-2}
E_1^{0,2} = H^2_{G \rel}(\bF), \;\;
E_1^{1,2} = H^1_{G \rel}(P_1),
\end{equation}
where the groups are again one-dimensional. $E_1^{0,2}$ is generated by $t$, and as a special case of the graded commutativity of the (relative) product, we have $d_1(t) = 0$. This immediately yields (ii). We now know that \eqref{eq:q-zero-0b} is nonzero in degree $2$; therefore, so is \eqref{eq:q-zero-0}, which provides the missing piece of (i).
\end{proof}

Let's make the consequence of Lemma \ref{th:q-zero-0}(ii) explicit. There is a cocycle of degree $2$,
\begin{equation} \label{eq:tq}
\begin{aligned}
& t_q \in Q^{h\Pi \rel}_0, \\
& t_q = t_0 + qt_1 + \cdots, \;\; t_m \in (P^{h\Pi\rel}_m)^{2-2m}, \\
& [t_0] \in H_{G \rel}^2(P_0) \iso H^2_G(\bF) \text{ is the generator $t$;}
\end{aligned}
\end{equation}
and the class $[t_q] \in H^*(Q^{h\Pi \rel}_0)$ is unique. By part (i), the same would be true if we allowed a more general choice of coefficients in $P^{h\Pi}$, but that will be less relevant for our future constructions.

\begin{remark} \label{th:shortcoming-1}
Our construction of $t_q$ depended on the fact that $t\lambda_G = \lambda_G t \in H^*_G(P_1)$. This is a shadow of a deeper homotopy commutativity, meaning that $t$ is homotopically a central element in $P^{h\Pi}$. Using that, one should be able to give a more elegant construction of the higher order terms in $t_q$, avoiding any kind of cohomological vanishing considerations. 

In discussions involving $G$-equivariant Floer cohomology, these kinds of problems are often bypassed by using the minimal model $\bF[t]$ instead of $\bF^{hG}$, since that is a strictly commutative algebra. Other $G$-equivariant structures are then constructed so as to be $\bF[t]$-linear. However, that setup is less convenient for comparing the theory with its counterpart for the finite group $\Gamma$. (It is instructive to look at \cite{sen23}, where the comparison involves passing to a larger cellular model for $G$ \cite[Definition 2.3.1]{sen23}, and hence for equivariant cochains; $t$-linearity is lost in that enlargement.)
\end{remark}

\begin{lemma} \label{th:q-zero-1}
(i) The map 
\begin{equation} \label{eq:q-zero-1a}
H^*(Q^{h\Pi}_p) \longrightarrow H^*_G(P_p) \iso H_{-*}(S^1)[t]
\end{equation}
(see \eqref{eq:p-splitting} for the second isomorphism) is an isomorphism in degrees $\leq 0$, as well as in degree $2$.

(ii) The map
\begin{equation} \label{eq:q-zero-1b}
H^*(Q_p^{h\Pi\rel}) \longrightarrow H^*_{G\rel}(P_p) \iso t H_{-*}(S^1)[t]
\end{equation}
is an isomorphism in degrees $\leq 2$.
\end{lemma}

\begin{proof}
The proof follows that of Lemma \ref{th:q-zero-0}. For (i), the nontrivial terms on the $E_1$ page of degree $\leq 3$ are (each one-dimensional)
\begin{equation} \label{eq:low-degree-3}
\begin{aligned}
& E_1^{0,-1} = H_G^{-1}(P_p), \;\; E_1^{0,0} = H_G^0(P_p), \;\;
E_1^{0,1} = H_G^1(P_p), \;\; E_1^{1,0} = H_G^{-1}(P_{p+1}), \;\;
\\
& E_1^{0,2} = H_G^2(P_p), \;\;
E_1^{0,3} = H_G^3(P_p), \;\; E_1^{2,1} = H_G^{-1}(P_{p+2}).
\end{aligned}
\end{equation}
One of the relevant differentials can be analyzed by looking at the diagram
\begin{equation} \label{eq:d1-diagram}
\xymatrix{
\ar[d]
E_1^{0,0} = H^0_G(P_p) \ar[rr]^-{d_1} && E_1^{1,0} = H^{-1}_G(P_{p+1}) \ar[d]
\\
H^0(P_p) \ar[rr]^-{x \mapsto \lambda x - x \lambda} && H^{-1}(P_{p+1})
}
\end{equation}
Here, the right $\downarrow$ is an isomorphism by Lemma \ref{th:equivariant-cohomology}(ii). The bottom $\rightarrow$ is zero by \eqref{eq:iota-lambda}. Hence, the $d_1$ in \eqref{eq:d1-diagram} vanishes. This shows that \eqref{eq:q-zero-1a} is an isomorphism in degrees $\leq 0$, and injective in degree $2$.

(ii) For the relative version, the corresponding groups (again, one-dimensional) reduce to
\begin{equation} \label{eq:low-degree-4}
E_1^{0,1} = H_{G \rel}^1(P_p), E_1^{0,2} = H_{G \rel}^2(P_p), \;\;
E_1^{0,3} = H_{G \rel}^3(P_p), \;\; E_1^{1,2} = H_{G \rel}^1(P_{p+1}).
\end{equation}
The differential $d_1: E_1^{0,2} \rightarrow E_1^{1,2}$ is obtained from that in \eqref{eq:d1-diagram} by applying the isomorphism from Lemma \ref{th:relative-t}(i), hence again vanishes.
This shows (ii), and retroactively also adds the missing piece of (i).
\end{proof}

As before, let's note some concrete consequences. We can choose a cocycle of degree $0$,
\begin{equation}
\begin{aligned}
\label{eq:rq}
& r_q \in Q^{h\Pi}_p, \\
& r_q = r_p + q r_{p+1} + \cdots, \;\; r_m \in P_m^{h\Pi,2p-2m}, \\
& [r_p] = \iota^p_G;
\end{aligned}
\end{equation} 
and another one of degree $-1$,
\begin{equation}
\label{eq:rhoq}
\begin{aligned}
& \rho_q \in Q^{h\Pi}_p, \\
& \rho_q = \rho_p + q\rho_{p+1} + \cdots, \;\; \rho_m \in P_m^{h\Pi,2p-2m-1}, \\
& [\rho_p] = (\iota^{p-1}\lambda)_G.
\end{aligned}
\end{equation}
The cohomology classes $[r_q] \in H^0(Q_p^{h\Pi})$, $[\rho_q] \in H^{-1}(Q_p^{h\Pi})$ are uniquely determined by these assumptions. These classes will play the role of quantum Steenrod operations in our abstract context. It is convenient to know that they commute with \eqref{eq:tq}:

\begin{lemma} \label{th:steenrod-is-t-linear}
The cocycles 
\begin{equation}
t_q r_q - r_q t_q, \;\; t_q \rho_q - \rho_q t_q \in Q^{h\Pi\rel}_p
\end{equation}
are nullhomologous.
\end{lemma}

\begin{proof}
The $q = 0$ term of $t_q r_q - r_q t_q$ represents $t \iota^p_G - \iota^p_G t \in H^2_{G\rel}(P_p)$, which is zero. From that, the desired result follows by applying Lemma \ref{th:q-zero-1}(ii). The second case is parallel.
\end{proof}

While \eqref{eq:rq} seems a natural candidate for a quantum Steenrod operation, as its leading order term is geometrically represented by a single cylinder with $p$ marked points, the appearance of \eqref{eq:rhoq} may be unexpected. We will spend a little more time explaining how it fits. Write $\tilde{P}^{h\Pi}$ for the analogue of $P^{h\Pi}$ defined by taking coefficients in $\tilde{\bF}$, see \eqref{eq:witt-vectors}. Lemma \ref{th:deltaq} goes through in this context, so that one can define an element $\tilde{\delta}_q$, and a deformed differential on the corresponding space $\tilde{Q}^{h\Pi}$. Supposing that $\delta_q$ is taken to be the reduction of $\tilde{\delta}_q$ to $\bF$, one gets a Bockstein operation, which fits into a commutative diagram with that from \eqref{eq:bockstein},
\begin{equation} \label{eq:bockstein-2}
\xymatrix{
\ar[d]
H^*(Q_p^{h\Pi}) \ar[rr]^-{\beta} && H^{*+1}(Q_p^{h\Pi}) \ar[d]
\\
H^*_G(P_p) \ar[rr]^{\beta} && H^{*+1}_G(P_p)
}
\end{equation}
One then has

\begin{lemma} \label{th:bockstein-3}
Under the Bockstein from \eqref{eq:bockstein-2}, one has
\begin{equation} \label{eq:bockstein-3}
[t_q \rho_q] = \beta([r_q]).
\end{equation}
\end{lemma}

\begin{proof}
From \eqref{eq:bockstein-2} and Lemma \ref{th:bockstein}, one sees that the difference between the two sides of \eqref{eq:bockstein-3} maps to zero in $H^1_G(P_p)$. This degree unfortunately falls right outside the scope of Lemma \ref{th:q-zero-1}(i); but we can work around that, as follows.

The left hand side of \eqref{eq:bockstein-3} obviously lifts to $H^*(Q^{h\Pi \rel}_p)$. On the right hand side, if we map $[r_q]$ to $H^*(Q_p)$, then that image admits a lift to $H^*(\tilde{Q}_p)$, by first lifting the leading order term $\iota^p$, and then carrying out one of our usual spectral sequence (or obstruction theory) arguments. As a consequence of that, $\beta([r_q])$ also maps to zero in $H^*(Q_p)$, which means that the class can be lifted to $H^*(Q^{h\Pi\rel}_p)$ as well. 

With that in mind, let's reconsider the difference between the two sides of \eqref{eq:bockstein-3}.
We know that it can be represented by a class in the kernel of the composition
\begin{equation}
H^1(Q^{h\Pi \rel}_p) \longrightarrow H^1_{G \rel}(P_p) \longrightarrow H^1_G(P_p).
\end{equation}
As a consequence of Lemma \ref{th:equivariant-cohomology}(iii), the second map is injective, so our class lies in the kernel of the first map. By Lemma \ref{th:q-zero-1}(ii), the class must be zero.
\end{proof}

\section{The connection\label{sec:connection}}
We begin by recalling the definition of $p$-curvature, in an abstract Weyl algebra framework. Then, returning to the framework from the previous section, we construct an abstract connection differentiating with respect to $q$, as an element in a suitable Weyl algebra \eqref{eq:universal-connection-1}. 

After that, we compute the $p$-curvature of that connection, in Theorem \ref{th:p-curvature-conjecture}. This may seem to be overly abstract, since all computations are carried out on complexes of operations, without anything that they act on. In the last part of this section, we therefore introduce a notion of module (an abstract counterpart of a Floer complex), and see what the implications for that object are.

\subsection{Weyl algebra background\label{subsec:weyl}}
Start with the Weyl algebra $\bF[[q]]\langle \dq \rangle$, where as usual $\bF = \bF_p$. We write $\dq$ to distinguish this element from the usual operation of $q$-differentiation of functions, denoted by $\partial_q$. The relation between the two is that
\begin{equation} \label{eq:weyl-relation}
\dq\cdot f - f\cdot \dq = \partial_q f \;\; \text{ for $f \in \bF[[q]]$}.
\end{equation}
Because we are working in characteristic $p$, the element $\dq^p$ is central. We divide by the ideal generated by it, forming the truncated Weyl algebra $\bF[[q]]\langle \dq \rangle/(\dq^p)$.
An element in that algebra can be written uniquely as 
\begin{equation} \label{eq:weyl-element}
x = \sum_{j=0}^{p-1} x_j \dq^j \quad\text{ with $x_0,\dots,x_{p-1} \in \bF[[q]]$.}
\end{equation}
Let's say that $x$ is a differential operator of order $\leq k$ if only the terms with $j \leq k$ in \eqref{eq:weyl-element} are nonzero. One can reformulate the definition more intrinsically as follows: $x$ is a differential operator of order $0$ if it commutes with $q$; and for $k>0$, $x$ is a differential operator of order $\leq k$ if $x q - q x$ is of order $\leq k-1$.

Given an associative algebra $A$ over $\bF$, one can similarly form the truncated Weyl algebra $A[[q]]\langle\dq \rangle/(\dq^p)$. Elements can be written in analogy with \eqref{eq:weyl-element}, where the coefficients now lie in $A[[q]]$, and the same notion of order carries over.
Given some $\alpha_q \in A[[q]]$, one writes an abstract connection
\begin{equation} \label{eq:connection}
\nablaq = \dq + \alpha_q \in A[[q]]\langle \dq \rangle/(\dq^p),
\end{equation}
which is a differential operator of order $1$. From the general formula
\begin{equation} \label{eq:k-th-order}
\nablaq^k q - q \nablaq^k = k\nablaq^{k-1}
\end{equation}
one sees that the $p$-th power $\nablaq^p$ is a differential operator of order $0$, hence belongs to $A[[q]]$. This is called the $p$-curvature of $\nablaq$. One can write an explicit formula for it, in terms of $\alpha_q$ and its $q$-derivatives (see \cite[Theorem 2]{schimming-rida96} or \cite[Section 3]{osserman08}), but we will not need that.

\subsection{The definition}
We return to the setup from Section \ref{sec:tqft}. One can apply term-by-term differentiation in the variable $q$ to elements of $Q^{h\Pi}$, which gives a map 
\begin{equation} \label{eq:differentiation}
\partial_q: Q^{h\Pi,j}_i \longrightarrow Q^{h\Pi,j-2}_{i+1}.
\end{equation}
In particular, by differentiating \eqref{eq:deltaq} one gets a cocycle $\partial_q \delta_q \in Q^{h\Pi,-1}_1$:
\begin{equation} \label{eq:derivative-of-delta}
\begin{aligned}
& \partial_q \delta_q = \delta_1 + 2q\delta_2 + \cdots, \\
& d_q(\partial_q\delta_q) = d(\partial_q\delta_q) + \delta_q (\partial_q \delta_q)
+ (\partial_q \delta_q) \delta_q = \partial_q (d\delta_q + \delta_q^2) = 0.
\end{aligned}
\end{equation}
This measures the failure of \eqref{eq:differentiation} to be a chain map. Namely, differentiating \eqref{eq:dq} yields
\begin{equation} \label{eq:derivative-of-differential}
\partial_q (d_q x) - d_q(\partial_q x) = (\partial_q\delta_q) x - (-1)^{|x|} x (\partial_q \delta_q)
\;\; \text{for any $x \in Q^{h\Pi}$.}
\end{equation}

\begin{lemma} \label{th:universal-connection}
(i) There is a degree $0$ cochain
\begin{equation} \label{eq:universal-connection}
\begin{aligned}
& \alpha_q \in Q^{h\Pi}_1, \\
& \alpha_q = \alpha_1 + q\alpha_2 + \cdots, \;\; \alpha_m \in P^{h\Pi,2-2m}_m, \\
& d_q\alpha_q = t_q (\partial_q \delta_q).
\end{aligned}
\end{equation}

(ii) For any such choice of $\alpha_q$, the image of $\alpha_1$ in $P_1$ represents $\iota$.
\end{lemma}

\begin{proof}
We will start by constructing only the first order part of (i). The equation is
\begin{equation}
d\alpha_1 = t_0 \delta_1.
\end{equation}
The right hand side represents $t \lambda_G \in H^1_G(P_1)$. We know from Lemma \ref{th:equivariant-cohomology}(ii) that this is zero, so one can definitely find some $\alpha_1$ with that property. 

Let's turn to (ii). By construction, $t_0$ maps to zero in $P_0$. Therefore, the image of $\alpha_1$ in $P_1$ is a cocycle, whose cohomology class maps to $t\lambda_G$ under the connecting map from \eqref{eq:relative-les}. As a consequence, if we consider \eqref{eq:three-maps} for $V = P_1$, meaning the composition
\begin{equation} \label{eq:three-maps-2}
H^0(P_1) \longrightarrow H^1_{G \rel}(P_1) \stackrel{t}{\longleftarrow} H^{-1}_G(P_1) \longrightarrow H^{-1}(P_1),
\end{equation}
then that maps the image of $\alpha_1$ to $\lambda$. By Lemma \ref{th:relative-t}(ii), the same holds for $\iota$, since loop rotation takes it to $\lambda$. All three maps in \eqref{eq:three-maps-2} being isomorphisms, the desired statement follows.

Now we return to the higher order parts of (i). The second order equation is
\begin{equation}
d\alpha_2 = t_1 \delta_1 + 2t_0 \delta_2 - \delta_1 \alpha_1 + \alpha_1 \delta_1,
\end{equation}
where the right hand side is a degree $-1$ cocycle in $P_2^{h\Pi}$. We know from Lemma \ref{th:equivariant-cohomology}(ii) that $H^{-1}_G(P_2) \rightarrow H^{-1}(P_2)$ is an isomorphism, hence can apply that forgetful map to our cocycle: by (ii), the outcome represents $\iota\lambda-\lambda\iota$, which is zero by Lemma \ref{th:cohomology}. The obstructions of order $m \geq 3$ lie in $H^{3-2m}(P_m^{h\Pi})$, and all those groups are zero.
\end{proof}

Next, we consider the compatibility between $\alpha_q$ and \eqref{eq:tq}. As a consequence of \eqref{eq:derivative-of-differential} and \eqref{eq:universal-connection}, the following is a degree $2$ cocycle:
\begin{equation} \label{eq:connection-is-t-linear}
t_q (\partial_qt_q) + \alpha_q t_q - t_q \alpha_q \in Q^{h\Pi \rel}_1.
\end{equation}

\begin{lemma} \label{th:connection-is-t-linear}
The cocycle \eqref{eq:connection-is-t-linear} is nullhomologous.
\end{lemma}

\begin{proof}
By definition, \eqref{eq:connection-is-t-linear} lies in $Q^{h\Pi \rel}_1$. We have not specifically considered the cohomology of that complex before, but a $q$-filtration argument together with Lemma \ref{th:equivariant-cohomology}(iv) shows that it vanishes in degree $2$.
\end{proof}

We now introduce differential operators into our context. For given $i,j \in \bZ$, let $W_i^{h\Pi,j}$ be the space of series
\begin{equation}
x_q = \sum_{\substack{0 \leq k \leq p-1 \\ m \geq i-k}} x_{m,k} q^{m-i+k} \dq^k, \;\; 
x_{m,k} \in P_m^{h\Pi,j-2m+2i}.
\end{equation}
One equips the bigraded space $W^{h\Pi}$ with a multiplicative structure, by combining that on $P^{h\Pi}$ with the truncated Weyl algebra. Note that 
\begin{equation}
\begin{aligned}
& qW_i^{h\Pi,j} \subset W_{i-1}^{h\Pi,j+2}, \\
& \dq W_i^{h\Pi,j} \subset W_{i+1}^{h\Pi,j-2}.
\end{aligned}
\end{equation}
We introduce a deformed differential on $W^{h\Pi}$ as in \eqref{eq:dq}. For instance,
writing $\dq = e \, \dq \in W^{h\Pi}_1$ where $e$ is the unit, one has
\begin{equation} \label{eq:dq-of-dq}
d_q(\dq) = \delta_q \dq - \dq \delta_q = -\partial_q \delta_q.
\end{equation}
The previous $Q^{h\Pi}$ is included in $W^{h\Pi}$ as a dg subalgebra, consisting of those elements that commute with $q$ (order zero differential operators). There is also a non-equivariant version $W$, and one defines $W^{h\Pi \rel}$ as the kernel of the forgetful map $W^{h\Pi} \rightarrow W$.

\begin{lemma} \label{th:order-zero}
(i) Let $x \in W^{h\Pi}_m$, $m>0$, be a cocycle such that
\begin{equation} \label{eq:y-discrepancy}
xq - qx = d_qy,
\end{equation}
where $y \in W^{h\Pi}_{m-1}$ is a differential operator of order $\leq p-2$. Then there is some $z \in W^{h\Pi}_m$ such that $x-d_qz \in Q^{h\Pi}_m$.

(ii) If $y$ can be chosen to lie in $W^{h\Pi \rel}_{m-1}$, then  one can also take $z \in W^{h\Pi \rel}_m$.
\end{lemma}

\begin{proof}
(i) Let's write concretely
\begin{equation} \label{eq:y-discrepancy-2}
 y = \sum_{j=0}^k y_j \dq^j, \;\; y_j \in Q^{h\Pi}_{m-1-j},
\end{equation}
where $k \leq p-2$. Take $z_k = y_k \frac{\dq^{k+1}}{k+1}$. By \eqref{eq:dq-of-dq}, this satisfies 
\begin{equation}
d_q z_k = \textstyle (d_q y_k) \frac{\dq^{k+1}}{k+1} + \text{(differential operator of order $\leq k$),}
\end{equation}
and therefore
\begin{equation}
\begin{aligned}
(d_qz_k) q - q (d_qz_k) & = (d_q y_k) \dq^k + \text{(differential operator of order $<k$)}
\\ & = d_q(y_k \dq^k) + \text{(differential operator of order $<k$)} 
\\ & = d_qy + \text{(differential operator of order $<k$)} 
\\ & = xq - qx + \text{(differential operator of order $<k$).}
\end{aligned}
\end{equation}
By subtracting $d_qz_k$ from $x$, one gets a homologous cocycle for which the counterpart of \eqref{eq:y-discrepancy-2} has $k-1$ instead of $k$. After repeating that $k$ times, one gets a cocycle which commutes with $q$, which is the desired outcome.

(ii) Is clear from inspecting the proof.
\end{proof}

Let's introduce our formal connection:
\begin{equation} \label{eq:universal-connection-1}
\nablaq \stackrel{\mathrm{def}}{=} t_q\dq + \alpha_q \in W_1^{h\Pi}.
\end{equation}
Then, \eqref{eq:universal-connection} says that
\begin{equation}
d_q(\nablaq) = 0,
\end{equation}
and Lemma \ref{th:connection-is-t-linear} has the following interpretation:
\begin{equation} \label{eq:connection-is-t-linear-2}
\begin{aligned}
t_q \nablaq - \nablaq t_q = d_q y\;\; \text{for some } y \in Q_1^{h\Pi \rel}.
\end{aligned}
\end{equation}

\begin{lemma} \label{th:order-zero-2}
There is a $z \in W^{h\Pi \rel}_p$ such that $\nablaq^p - d_qz$ is a differential operator of order $0$ (meaning, it lies in $Q^{h\Pi}_p$).
\end{lemma}

\begin{proof}
This is a version of the elementary computation \eqref{eq:k-th-order}, but we have to be more careful, since $t_q$ is not strictly $\nablaq$-covariant. From the definition of $\nablaq$, one gets
\begin{equation}
\begin{aligned}
& \nablaq^k q - q \nablaq^k = \sum_{0 \leq j \leq k-1} \nablaq^j t_q \nablaq^{k-1-j} \\
& = k t_q \nablaq^{k-1} + \sum_{0 \leq j \leq k-2} (k-1-j) \nablaq^j (\nablaq t_q - t_q \nablaq) \nablaq^{k-2-j}.
\end{aligned}
\end{equation}
Specializing to $k = p$, and applying \eqref{eq:connection-is-t-linear-2}, yields
\begin{equation}
\nablaq^p q - q \nablaq^p  = d_q w, \quad w = -\sum_{0 \leq j \leq p-2} (p-1-j) \nablaq^j y \nablaq^{p-2-j}.
\end{equation}
The right hand side lies in $W^{h\Pi \rel}$, and is a differential operator of order $\leq p-2$. 
Lemma \ref{th:order-zero} then concludes the argument.
\end{proof}

\begin{theorem} \label{th:p-curvature-conjecture}
We have 
\begin{equation}
[\nablaq^p] =  [r_q] \in H^0(W_p^{h\Pi}).
\end{equation}
\end{theorem}

\begin{proof}
Consider the diagram
\begin{equation} \label{eq:forget-forget}
\xymatrix{
H^0(Q_p^{h\Pi}) \ar[rr]^{\iso} \ar[d] && H^0(P_p^{h\Pi}) \ar[d]^{\iso} \\
H^0(Q_p) \ar[rr] && H^0(P_p)
}
\end{equation}
where the isomorphisms are established by Lemmas \ref{th:q-zero-1} and \ref{th:equivariant-cohomology}(iii). Start with the cocycle $\nablaq^p - d_qz \in Q_p^{h\Pi}$ from Lemma \ref{th:order-zero-2}. By definition, the image of this under the forgetful map to $Q_p$ equals $\alpha_q^p$. If we then further map that to $P_p$, the outcome is $\alpha_1^p$, which by Lemma \ref{th:universal-connection}(ii) represents $\iota^p$. Now, if I take $r_q$ and apply the same forgetful maps, then one gets a representative of $\iota^p$, by definition \eqref{eq:rq}. This proves that $[\nablaq^p - d_qz] = [r_q] \in H^0(Q_p^{h\Pi})$. Mapping that equality to $W_p^{h\Pi}$ yields the desired result.
\end{proof}

\subsection{Application to dg modules\label{subsec:module}}
We investigate the following general situation:

\begin{assumptions} \label{th:module}
Let $E$ be a chain complex which comes with a $G$-action, in the sense of Section \ref{subsec:topological-groups}. We restrict that to the action of $\Pi \subset G$, and ask $E$ to also have a compatible structure of a module over \eqref{eq:p}, given by $\Pi$-equivariant chain maps
\begin{equation} \label{eq:module}
P_m \otimes E \longrightarrow E[2m].
\end{equation}
\end{assumptions}

As a consequence, $E^{h\Pi}$ becomes a module over $P^{h\Pi}$, by 
\begin{equation} \label{eq:equivariant-module}
P_m^{h\Pi} \otimes E^{h\Pi} \stackrel{\times}{\longrightarrow} (P_m \otimes E)^{h\Pi}
\stackrel{\eqref{eq:module}}{\longrightarrow} E^{h\Pi}[2m].
\end{equation}
Equip the space
\begin{equation} \label{eq:f-space}
F^{h\Pi} \stackrel{\mathrm{def}}{=} E^{h\Pi}[[q]]
\end{equation}
with the $q$-linear differential
\begin{equation} \label{eq:f-differential}
d_qx = dx + \delta_q x.
\end{equation}

\begin{remark}
The grading in \eqref{eq:module} is so that in combination with \eqref{eq:dq}, the differential \eqref{eq:f-differential} has degree $1$. There are other options, where the degree of \eqref{eq:module} scales differently with $m$, and one correspondingly assigns a nontrivial degree to the variable $q$. These correspond to different situations in symplectic topology (ours being the Calabi-Yau case).
\end{remark}

By the $q$-linear extension of \eqref{eq:equivariant-module}, $F^{h\Pi}$ becomes a dg module over $Q^{h\Pi}$; note that while $Q^{h\Pi}$ is bigraded, the same is not true of \eqref{eq:f-space}. In particular, $H^*(F^{h\Pi})$ becomes a graded $\bF[t]$-module by the action \eqref{eq:tq},
\begin{equation} \label{eq:tq-action}
[t_q]: H^*(F^{h\Pi}) \longrightarrow H^{*+2}(F^{h\Pi}).
\end{equation}
The action of \eqref{eq:rq} and \eqref{eq:rhoq} yields endomorphisms
\begin{align}
\label{eq:rq-action}
& [r_q]: H^*(F^{h\Pi}) \longrightarrow H^{*+2p}(F^{h\Pi}), \\
\label{eq:rhoq-action}
& [\rho_q]: H^*(F^{h\Pi}) \longrightarrow H^{*+2p-1}(F^{h\Pi}), 
& \end{align}
These operations are $\bF[[q]]$-linear (by definition) and $\bF[t]$-linear (by Lemma \ref{th:steenrod-is-t-linear}). 

For future purposes, we need one more observation. The situation here is that of Lemma \ref{th:bockstein-3}, which means that $\delta_q$ is the reduction of some $\tilde{\delta}_q$ defined over $\tilde{\bF}$. This gives rise to a Bockstein on $H^*(Q^{h\Pi})$.

\begin{assumptions} \label{th:module-2}
(i) The structure of $E$ (meaning, the $G$-action, and also \eqref{eq:module} which gives the $\Pi$-equivariant dg module structure) can be lifted from $\bF$ to $\tilde{\bF}$. By this, we mean that there is a chain complex of free $\tilde{\bF}$-modules $\tilde{E}$, which has the corresponding structures, and whose mod $p$ reduction recovers $E$. 

(ii) The lift from (i) is such that the induced Bockstein on $H^*(F^{h\Pi})$ is zero.

(iii) Moreover, $H^*(F^{h\Pi})$ is a torsion-free $\bF[t]$-module.
\end{assumptions}

\begin{lemma} \label{th:no-rho}
Under Assumptions \ref{th:module-2}, the endomorphism \eqref{eq:rhoq-action} is zero.
\end{lemma}

\begin{proof}
The Bockstein $\beta$ on $H^*(F^{h\Pi})$ is compatible with the action of $H^*(Q^{h\Pi})$, in the sense that
\begin{equation} \label{eq:bockstein-module}
\beta(xy) = \beta(x)y + (-1)^{|x|} x\beta(y) \quad \text{ for $x \in H^*(Q^{h\Pi})$, $y \in H^*(F^{h\Pi})$.}
\end{equation}
Under our assumption (ii), we have $\beta(y) = 0$ and $\beta(xy) = 0$, which means that $\beta(x)$ acts trivially on $H^*(F^{h\Pi})$. By Lemma \ref{th:bockstein-3}, this means that the composition of $t$ and \eqref{eq:rhoq-action} is zero. But by assumption (iii), the first of these maps is injective, so the second one must vanish.
\end{proof}

We make $F^{h\Pi}$ into a module over $W^{h\Pi}$, by declaring that $\dq$ acts by $q$-differentiation. The action of the cocycle $\nablaq$ from \eqref{eq:universal-connection-1} then defines a chain map
\begin{equation}
\begin{aligned}
& F^{h\Pi} \longrightarrow F^{h\Pi}[2], \\
& y \longmapsto \nablaq y = t_q \partial_q y + \alpha_q y.
\end{aligned}
\end{equation}
We write $\nabla_q$ for the induced action on $H^*(F^{h\Pi})$, which satisfies the Leibniz rule
\begin{equation}
\nabla_q (qx) = q\nabla_q x + t x.
\end{equation}
It it also $\bF[t]$-linear, as a consequence of \eqref{eq:connection-is-t-linear}. Finally, from Theorem \ref{th:p-curvature-conjecture} we immediately get:

\begin{corollary} \label{th:p-curvature-module}
For any module as in Assumptions \ref{th:module}, the $p$-th power of the connection $\nabla_q$ on $H^*(F^{h\Pi})$ equals \eqref{eq:rq-action}.
\end{corollary}

\section{Finite symmetries\label{sec:tqft-2}}
The usual theory of quantum Steenrod operations is based on a finite group of rotations of the cylinders (see e.g.\ \cite{seidel-wilkins21}). In this section, we set up a corresponding version of the formalism from Section \ref{sec:tqft}. The approach is parallel to the previous one (and we will reflect that by omitting proofs), but the outcome is different in a couple of respects. First, because the version with finite symmetry has a larger coefficient ring, there is in principle a larger space of choices; one instance of that is the clarification we have to make when choosing \eqref{eq:iota-gamma}. More importantly, the behaviour of the Bockstein is different: there is no counterpart of Lemma \ref{th:bockstein-3} for finite symmetry groups.

At the end of the section, we pick up the comparison between $S^1$-equivariance and that for finite groups. Under suitable assumptions, this leads us to a description of the $p$-curvature in the form suggested by the original (Etingof-Lee) conjecture.

\subsection{Parameter spaces revisited\label{subsec:revisited}}
As in Section \ref{subsec:finite-cyclic-group}, $\Gamma \subset S^1$ is the cyclic subgroup of (odd prime) order $p$. We will continue working over $\bF =\bF_p$ throughout. We use the same spaces $\bar\frakP_m$ and chain complexes $P_m$ as before. The equivariant version
\begin{equation}
P^{h\Gamma} = \bigoplus_{m \geq 0} P_m^{h\Gamma}
\end{equation}
is a bigraded dga, which comes with a forgetful map $P^{h\Gamma} \rightarrow P$. On the cohomology level, 
\begin{equation} \label{eq:p-forget-2}
H^*(P_m^{h\Gamma}) = H^*_{\Gamma}(P_m) \longrightarrow H^*(P_m).
\end{equation}

\begin{lemma} \label{th:equivariant-cohomology-2}
(i) $H^*_\Gamma(P_0) \iso H^*_{\Gamma}(\bF) = \bF[t,\theta]$.

(ii) If $m$ is not a multiple of $p$, we have $H^*_\Gamma(P_m) \iso H_{-*}(S^1)$, and
\begin{equation} \label{eq:pm-theta}
\theta: H^{-1}_{\Gamma}(P_m) \stackrel{\iso}{\longrightarrow} H^0_{\Gamma}(P_m).
\end{equation}
The map \eqref{eq:p-forget-2} is an isomorphism in degrees $\leq -1$, but is zero in degree $0$.

(iii) If $m>0$ is a multiple of $p$, then $H^*_{\Gamma}(P_m)$ is a free graded module over $\bF[t,\theta]$, with one generator each in degrees $-1$ and $0$. The map \eqref{eq:p-forget-2} is an isomorphism in degrees $\leq -1$, and surjective in degree $0$.

(iv) For any $m$, $H^*_{\Gamma \rel}(P_m)$ is concentrated in degrees $\geq 0$; and if $m$ is not a multiple of $p$, it is concentrated in degrees $\{0,1\}$.
\end{lemma}

One could derive this from its $S^1$-equivariant counterpart (Lemma \ref{th:equivariant-cohomology}) by applying the general Lemma \ref{th:comparison}, but we prefer to give a direct proof; besides being more geometrically intuitive, we will refer to one piece of the argument later on.

\begin{proof}
(i) is obvious. So is (iii) because then, the $\Gamma$-action on $\bar\frakP_m$ corresponds to the trivial action on the target $S^1$ in \eqref{eq:morton}.

In the remaining case (ii), the $\Gamma$-action on $\bar\frakP_m$ is free, so that $H^*_{\Gamma}(P_m) \iso H^*(P_m^{\Gamma}) = H_{-*}(\bar\frakP_m/\Gamma)$. In that context, the map \eqref{eq:pm-theta} is the cap action of the class in $H^1(\bar\frakP_m/\Gamma)$ which classifies the covering $\bar\frakP_m \rightarrow \bar\frakP_m/\Gamma$. Under \eqref{eq:morton}, this corresponds to a nontrivial $p$-fold cover of $S^1$, which explains \eqref{eq:pm-theta}. In the same terms, the map \eqref{eq:p-forget-2} is the transfer $H_{-*}(\bar\frakP_m/\Gamma) \rightarrow H_{-*}(\bar\frakP_m)$, hence obviously has the properties we have stated.
\end{proof}

As before, we write $(\lambda\iota^{m-1})_{\Gamma}$ for the unique equivariant lift of $\lambda\iota^{m-1}$. In case (iii) of Lemma \ref{th:equivariant-cohomology-2}, the choice of lift of $\iota^m$ to a class
\begin{equation} \label{eq:iota-gamma}
\iota^m_{\Gamma} \in H^0_{\Gamma}(P_m)
\end{equation}
is not unique. We choose the lift represented by the element of $P_m^{\Gamma} \subset P_m^{h\Gamma}$ arising from the invariant configuration
\begin{equation} \label{eq:symmetric-configuration}
(\overbrace{(0,0), \dots, (0,0)}^{\text{$m/p$ times}}, \overbrace{(0,1/p),\dots,(0,1/p)}^{\text{$m/p$ times}},\dots) \in \mathit{Sym}_m(\bR \times S^1)^{\Gamma}.
\end{equation}
These lifts are compatible with the product,
\begin{equation}
\iota^m_{\Gamma} \iota^n_{\Gamma} = \iota^{m+n}_{\Gamma} \in H^0_{\Gamma}(P_{m+n}).
\end{equation}

\subsection{The deformed complex}
The counterpart of Lemma \ref{th:deltaq} holds, for the same reason:

\begin{lemma} \label{th:deltaq-2}
(i) There is a formal series
\begin{equation} \label{eq:deltaq-2}
\begin{aligned}
& \delta_q = \delta_1 q + \delta_2 q^2 + \cdots, \;\; \delta_i \in P_i^{h\Gamma,1-2i}, \\
& d\delta_q + \delta^2_q = 0, \\
& [\delta_1] = \lambda_{\Gamma} \in H^{-1}_{\Gamma}(P_1).
\end{aligned}
\end{equation}

(ii) Given $\delta_q$ and $\delta_q'$ both having the properties above, there is
\begin{equation} \label{eq:phiq-2}
\begin{aligned}
& \phi_q = \phi_1 q + \phi_2 q^2 + \cdots, \;\; \phi_m \in P_m^{h\Gamma,-2m}, \\
& d\phi_q + \delta_q' (e+\phi_q) - (e+\phi_q) \delta_q = 0. 
\end{aligned}
\end{equation}
\end{lemma}

We use that to define $Q^{h\Gamma}$ along the same lines as before, as well as the relative version $Q^{h\Gamma \rel}$. Next come the counterparts of Lemmas \ref{th:q-zero-0} and \ref{th:q-zero-1}, with proofs omitted.

\begin{lemma} \label{th:q-zero-2}
(i) The map (obtained by setting $q = 0$)
\begin{equation} \label{eq:q-zero-2}
H^*(Q^{h\Gamma}_0) \longrightarrow H^*_{\Gamma}(P_0) \iso \bF[t,\theta]
\end{equation}
is an isomorphism in degrees $\leq 0$, and surjective in degree $1$.

(ii) Similarly, 
\begin{equation} \label{eq:q-zero-2b}
H^*(Q^{h\Gamma \rel}_0) \longrightarrow H^*_{\Gamma \rel}(P_0) \iso \theta\bF[t] \oplus t\bF[t]
\end{equation}
is an isomorphism in degrees $\leq 1$.
\end{lemma}

\begin{lemma} \label{th:q-zero-3}
(i) The map
\begin{equation}
H^*(Q^{h\Gamma}_p) \longrightarrow 
H^*_{\Gamma}(P_p) \iso H_{-*}(S^1)[t,\theta]
\end{equation}
is an isomorphism in degrees $\leq 0$.

(ii) The map
\begin{equation}
H^*(Q^{h\Gamma \rel}_p) \longrightarrow 
H^*_{\Gamma \rel}(P_p) \iso \theta H_{-*}(S^1)[t] \oplus t H_{-*}(S^1)[t]
\end{equation}
is an isomorphism in degrees $\leq 1$.
\end{lemma}

Given Lemma \ref{th:q-zero-2}(ii), there is a degree $1$ cocycle
\begin{equation} \label{eq:thetaq}
\begin{aligned}
& \theta_q \in Q^{h\Gamma \rel}_0, \\
& 
\theta_q = \theta_0 + q\theta_1 + \cdots, \;\; \theta_m \in (P_m^{h\Gamma \rel})^{1-2m}, \\
&
[\theta_0] \in H_{\Gamma \rel}^1(P_0) = H_{\Gamma}^1(\bF) \text{ is the generator $\theta$;}
\end{aligned}
\end{equation}
and its class $[\theta_q] \in H^*(Q_0^{h\Gamma \rel})$ is uniquely determined by those conditions. Similarly, Lemma \ref{th:q-zero-3}(i) provides us with cocycles of degrees $0$ and $-1$,
\begin{equation} \label{eq:sq}
\begin{aligned}
& s_q \in Q^{h\Gamma}_p, \\
& s_q = s_p + q s_{p+1} + \cdots, \;\; s_m \in P_m^{h\Gamma,2p-2m}, \\
& [s_p] = \iota^p_{\Gamma},   
\end{aligned}
\end{equation} 
respectively
\begin{equation}
\begin{aligned} \label{eq:sigmaq}
& \sigma_q \in Q^{h\Gamma}_p, \\
& \sigma_q = \sigma_p + q\sigma_{p+1} + \cdots, \;\; \sigma_m \in P_m^{h\Gamma,2p-1-2m}, \\
& [\sigma_p] = (\lambda\iota^{p-1})_{\Gamma}.
\end{aligned}
\end{equation}
Lemma \ref{th:q-zero-3}(ii) can then be used to prove the analogue of Lemma \ref{th:steenrod-is-t-linear}:

\begin{lemma} \label{th:steenrod-is-theta-linear}
The cocycles 
\begin{equation}
\theta_q s_q - s_q \theta_q, \;\; \theta_q \sigma_q + \sigma_q \theta_q \in Q^{h\Gamma\rel}_p
\end{equation}
are nullhomologous.
\end{lemma}

\begin{remark} \label{th:shortcoming-2}
The reader will have noticed that we have lifted the action of $\theta$ to $Q^{h\Gamma}$, but have not done the same for $t$. Lemma \ref{th:q-zero-2} is not suitable for that, because the desired argument lies outside the range of degrees allowed there. One can introduce a different relative version of $\Gamma$-equivariant cohomology, which not only allows one to define an action of $t$, but also to introduce a connection $\nabla_q$ and to compute its $p$-curvature. We will now give a short account of that approach. It represents a path not taken in the body of the paper, but could be convenient for applications to $p$-th iterates of symplectic automorphisms, where no larger symmetry is present.

Let $V$ be a complex with a $\Gamma$-action. Pick a generator $\gamma \in \Gamma$, and write $H^*_{\gamma}(V)$ for equivariant cohomology with respect to the group $\bZ$ acting by iterates of $\gamma$. Instead of the full bar complex, we represent this on the chain level by
\begin{equation} \label{eq:truncated-complex}
V^{h\gamma} \stackrel{\mathrm{def}}{=} V \oplus V[-1], 
\end{equation}
with differential
\begin{equation} \label{eq:truncated-differential}
(v,\bar{v}) \longmapsto (dv, -d\bar{v} + \gamma v - v)
\end{equation}
($d$ being the differential on $V$). This has an exterior product analogous to \eqref{eq:chain-exterior-product},
\begin{equation} \label{eq:truncated-product}
(v,\bar{v}) \times (w,\bar{w}) =  (v \otimes w, (-1)^{|v|} v \otimes \bar{w} + \bar{v} \otimes \gamma w).
\end{equation}
The pullback map $H^*_{\Gamma}(V) \rightarrow H^*_{\gamma}(V)$ is realized by
\begin{equation} \label{eq:truncation}
V^{h\Gamma} \longrightarrow V^{h\gamma}, \quad
\phi \longmapsto (\phi^0, \phi^1(\gamma));
\end{equation}
which is surjective and compatible with products. It maps our cocycle representative for $\theta$ to $(0,1)$ (so that one could also write $V^{h\gamma} = V[\theta]$); and the representative for $t$ goes to zero. Write $V^{h\Gamma \rel \gamma}$ for the kernel of \eqref{eq:truncation}, and $H^*_{\Gamma \rel \gamma}(V)$ for its cohomology. In the case of trivial coefficients,
\begin{equation}
H^*_{\gamma}(\bF) = \bF[\theta], \quad H^*_{\Gamma \rel \gamma}(\bF) = t\bF[t,\theta].
\end{equation}
As one can see from the second of these expressions, the new relative theory behaves quite similarly to \eqref{eq:relative-g}. On the chain level, one gets an ideal $Q^{h\Gamma \rel \gamma} \subset Q^{h\Gamma}$ which naturally accommodates a $\Gamma$-counterpart of \eqref{eq:tq}.
\end{remark}

The application to dg modules is strictly parallel to Section \ref{subsec:module}. We go through it only in order to fix notation. Let $E$ be a chain complex which comes with a module structure \eqref{eq:module}, equivariant with respect to $\Gamma$. Then, $E^{h\Gamma}$ becomes a module over $P^{h\Gamma}$ (with the same grading shifts as before). Similarly, 
\begin{equation} \label{eq:f-space-2}
F^{h\Gamma} \stackrel{\mathrm{def}}{=} E^{h\Gamma}[[q]]
\end{equation}
equipped with a suitable differential, becomes a dg module over $Q^{h\Gamma}$. In particular,  the action \eqref{eq:thetaq} gives an endomorphism
\begin{equation} \label{eq:thetaq-action}
[\theta_q]: H^*(F^{h\Gamma}) \longrightarrow H^{*+1}(F^{h\Gamma}).
\end{equation}
Similarly, \eqref{eq:sq} and \eqref{eq:sigmaq} yield $\theta$-linear endomorphisms
\begin{align}
& [s_q]: H^*(F^{h\Gamma}) \longrightarrow H^{*+2p}(F^{h\Gamma}), \\
\label{eq:sigmaq-action}
& [\sigma_q]: H^*(F^{h\Gamma}) \longrightarrow H^{*+2p-1}(F^{h\Gamma}).
\end{align}

\subsection{Comparison of equivariant theories}
We continue working over a field $\bF$ of odd characteristic $p$. The inclusion $\Gamma \subset \Pi$ gives rise to a homomorphism (of bigraded dga's)
\begin{equation} \label{eq:compare-p}
P^{h\Pi} \longrightarrow P^{h\Gamma}.
\end{equation}
We need to explain how it relates the generators from Section \ref{subsec:parameter-spaces} to their counterparts in Section \ref{subsec:revisited}.

\begin{lemma} \label{th:compare-iota-lambda}
Under the cohomology level map induced by \eqref{eq:compare-p},
\begin{align}
\label{eq:g-gamma-lambda}
& (\iota^{p-1}\lambda)_G \longmapsto (\iota^{p-1}\lambda)_{\Gamma}, \\
\label{eq:g-gamma-iota}
& \iota^p_G \longmapsto \iota^p_{\Gamma} + \theta (\iota^{p-1}\lambda)_{\Gamma}.
\end{align}
\end{lemma}

\begin{proof}
This is obvious for \eqref{eq:g-gamma-lambda}, since the two sides are the unique lifts of $\iota^{p-1}\lambda \in H^{-1}(P_p)$ to the respective equivariant cohomology groups.

For \eqref{eq:g-gamma-iota}, we know that $\iota^p_G$ is mapped to $\iota^p_{\Gamma} + c \theta (\iota^{p-1}\lambda)_{\Gamma}$ for some $c \in \bF$. Now apply the Bockstein:
\begin{align}
\label{eq:bockstein-4}
& \beta(\iota^p_G) = t (\lambda \iota^{p-1})_G, \\
\label{eq:bockstein-5}
& \beta(\iota^p_{\Gamma} + c \theta (\iota^{p-1}\lambda)_{\Gamma}) = c t \,(\iota^{p-1}\lambda)_{\Gamma}.
\end{align}
Here, \eqref{eq:bockstein-4} is Lemma \ref{th:bockstein}. Recall that there is a homotopy equivalence $\frakP_p \rightarrow S^1$, where the right hand side carries the trivial $\Gamma$-action. Hence,
\begin{equation} \label{eq:gamma-splitting}
H^*_{\Gamma}(P_p) \iso H^*_{\Gamma}(\bF) \otimes H_{-*}(S^1)
\end{equation}
Under this isomorphism, $\iota^p_{\Gamma}$ corresponds to $1 \otimes [\mathit{point}]$, as one can see by looking at the representative \eqref{eq:symmetric-configuration}, and therefore $\beta(\iota^p_{\Gamma}) = 0$. Still looking at \eqref{eq:gamma-splitting}, $\theta (\iota^{p-1}\lambda)_{\Gamma}$ corresponds to $\theta \otimes [S^1]$, and therefore $\beta(\theta (\iota^{p-1}\lambda)_{\Gamma}) = t (\iota^{p-1}\lambda)_{\Gamma}$ by \eqref{eq:beta-theta}; which explains \eqref{eq:bockstein-5}.

The restriction map from $G$-equivariant to $\Gamma$-equivariant cohomology is obviously compatible with Bocksteins. By applying that restriction map to \eqref{eq:bockstein-4}, and comparing the outcome with \eqref{eq:bockstein-5}, one sees that $c = 1$.
\end{proof}

Fix an element $\delta_q$ as in Lemma \ref{th:deltaq}, and the resulting dga $Q^{h\Pi}$. Under restriction to $\Gamma$, this gives a corresponding element satisfying the conditions from \ref{th:deltaq-2}, and if we use that one, the $q$-linear extension of \eqref{eq:compare-p} is a map of dga's,
\begin{equation} \label{eq:compare-q}
Q^{h\Pi} \longrightarrow Q^{h\Gamma}.
\end{equation}

\begin{lemma} \label{th:compare-q}
The restriction map \eqref{eq:compare-q} and the action of $\theta_q$ induce an isomorphism
\begin{equation} \label{eq:compare-q-2}
H^*(Q^{h\Pi}_m) \oplus \theta H^*(Q^{h\Pi}_m) \longrightarrow H^*(Q^{h\Gamma}_m).
\end{equation}
\end{lemma}

\begin{proof}
This requires a moment's thought, because of Remark \ref{th:prufer}. The underlying chain map is compatible with $q$-filtrations. The induced map between the $E_1$ pages of the associated spectral sequences is given by $H^*_G(P_i) \oplus \theta H^*_G(P_i) \rightarrow H^*_{\Gamma}(P_i)$, for $i \geq m$. Lemma \ref{th:comparison} shows that these maps are isomorphisms, and the desired statement follows from that.
\end{proof}

\begin{lemma} \label{th:compare-operations}
Under \eqref{eq:compare-q-2}, we have
\begin{align}
& [r_q] \longmapsto [s_q + \theta_q\sigma_q], \\
& [\rho_q] \longmapsto [\sigma_q].
\end{align}
\end{lemma}

\begin{proof}
By Lemma \ref{th:q-zero-3}(i), the cohomology classes on the right are determined by their $q = 0$ terms. This reduces the statement to Lemma \ref{th:compare-iota-lambda}.
\end{proof}

Suppose that we are in the situation of Assumptions \ref{th:module}. Then, a similar comparison process applies: the constructions \eqref{eq:f-space} and \eqref{eq:f-space-2} are related by a map 
\begin{equation}
F^{h\Pi} \longrightarrow F^{h\Gamma},
\end{equation}
and that induces an isomorphism
\begin{equation} \label{eq:compare-f}
H^*(F^{h\Pi}) \oplus \theta H^*(F^{h\Pi}) \stackrel{\iso}{\longrightarrow} H^*(F^{h\Gamma}).
\end{equation}

\begin{corollary} \label{th:p-curvature-module-2}
Suppose that Assumptions \ref{th:module-2} are satisfied. Take the obvious extension of $\nabla_q$ to an endomorphism of $H^*(F^{h\Pi}) \oplus \theta H^*(F^{h\Pi})$, meaning that $\nabla_q (\theta x) = \theta (\nabla_qx)$. Then, its $p$-th iterate fits into a commutative diagram
\begin{equation}
\xymatrix{
\ar[d]_-{\nabla_q^p}
H^*(F^{h\Pi}) \oplus \theta H^*(F^{h\Pi}) \ar[rr]_-{\iso}^-{\eqref{eq:compare-f}} && H^*(F^{h\Gamma})
\ar[d]^-{[s_q]}
\\
H^{*+2p}(F^{h\Pi}) \oplus \theta H^{*+2p}(F^{h\Pi}) \ar[rr]_-{\iso}^-{\eqref{eq:compare-f}} && H^{*+2p}(F^{h\Gamma}).
}
\end{equation}
\end{corollary}

\begin{proof}
From Corollary \ref{th:p-curvature-module}, we know that the left $\downarrow$ agrees with the action of $[r_q]$. By Lemma \ref{th:no-rho}, this is the same as the action of $[r_q] - \theta [\rho_q]$. Lemma \ref{th:compare-operations} implies that under \eqref{eq:compare-f}, this corresponds to the action of $[s_q] + [\theta_q\sigma_q] - [\theta_q\sigma_q] = [s_q]$.
\end{proof}


\begin{thebibliography}{10}

\bibitem{abouzaid-groman-varolgunes22}
M.~Abouzaid, Y.~Groman, and U.~Varolgunes.
\newblock Framed {$E_2$} structures in {F}loer theory.
\newblock {\em Adv. Math.}, 450:Paper No. 109755, 2024.

\bibitem{bai-lee-pomerleano26}
S.~Bai, J.~Lee, and D.~Pomerleano.
\newblock Arithmetic geometry of quantum connections on {C}alabi-{Y}au
  $3$-folds.
\newblock Preprint arXiv:2601.01654.

\bibitem{chen26}
Z.~Chen.
\newblock The {G}etzler-{G}auss-{M}anin connection and {K}ontsevich-{S}oibelman
  operations on the periodic cyclic homology.
\newblock Preprint arXiv:2601.16437.

\bibitem{chen24b}
Z.~Chen.
\newblock On operadic open-closed maps in characteristic $p$.
\newblock Preprint arXiv:2402.06183.

\bibitem{edidin-graham97}
D.~Edidin and W.~Graham.
\newblock Equivariant intersection theory.
\newblock {\em Invent. Math.}, 131:595--634, 1998.

\bibitem{etingof-varchenko24}
P.~Etingof and A.~Varchenko.
\newblock $p$-curvature of periodic pencils of flat connections.
\newblock Preprint arXiv:2401.05652.

\bibitem{katz70}
N.~Katz.
\newblock Nilpotent connections and the monodromy theorem: {A}pplications of a
  result of {T}urrittin.
\newblock {\em Inst. Hautes \'Etudes Sci. Publ. Math.}, 39:175--232, 1970.

\bibitem{lee23b}
J.~H. Lee.
\newblock Quantum {S}teenrod operations of symplectic resolutions.
\newblock {\em Geom. Topol.}, 29:4341--4387, 2025.

\bibitem{maclane}
S.~MacLane.
\newblock {\em Homology}.
\newblock Springer, 1963.

\bibitem{morton67}
H.~Morton.
\newblock Symmetric products of the circle.
\newblock {\em Proc. Camb. Phil. Soc.}, 63:349--352, 1967.

\bibitem{nikolaus-scholze18}
T.~Nikolaus and P.~Scholze.
\newblock On topological cyclic homology.
\newblock {\em Acta Math.}, 221:203--409, 2018.

\bibitem{osserman08}
B.~Osserman.
\newblock Frobenius-unstable bundles and {$p$}-curvature.
\newblock {\em Trans. Amer. Math. Soc.}, 360:273--305, 2008.

\bibitem{pomerleano-seidel26}
D.~Pomerleano and P.~Seidel.
\newblock The quantum connection and its mod {$p$} reduction.
\newblock Preprint arXiv:2606.28256.

\bibitem{pomerleano-seidel23}
D.~Pomerleano and P.~Seidel.
\newblock The quantum connection, {F}ourier-{L}aplace transform, and families
  of {$A_\infty$}-categories.
\newblock Preprint arXiv:2308.13567.

\bibitem{pomerleano-seidel24}
D.~Pomerleano and P.~Seidel.
\newblock Symplectic cohomology relative to a smooth anticanonical divisor.
\newblock Preprint arXiv:2408.09039.

\bibitem{rezchikov26}
S.~Rezchikov.
\newblock Noncommutative {C}artier formulae.
\newblock Preprint arXiv:2607.05360.

\bibitem{schimming-rida96}
R.~Schimming and S.~Rida.
\newblock Noncommutative {B}ell polynomials.
\newblock {\em Internat. J. Algebra Comput.}, 6:635--644, 1996.

\bibitem{seidel-wilkins21}
P.~Seidel and N.~Wilkins.
\newblock Covariant constancy of quantum {S}teenrod operations.
\newblock {\em J. Fixed Point Theory Appl.}, 24:Paper No. 52, 2022.

\bibitem{sen23}
D.~Sen.
\newblock {\em The $\mathbb{Z}/p\mathbb{Z}$-{G}ysin sequence for symplectic
  cohomology}.
\newblock PhD thesis, Univ. of Southern California, 2023.

\end{thebibliography}

\end{document}